\documentclass[leqno,12pt]{article}
\usepackage{amsmath,amssymb,amsthm,latexsym}
\def\min{\mathop{\mathrm{min}}\nolimits}

\def\proj{\mathop{\mathrm{proj}}\nolimits}

\def\op{\mathop{\mathrm{op}} \nolimits}
\def\aff{\mathop{\mathrm{aff}} \nolimits}

\def\id{\mathop{\mathrm{id}} \nolimits}

\newtheorem{dehilf}{Definition}[section]
\newtheorem{lehilf}[dehilf]{Lemma}
\newtheorem{sahilf}[dehilf]{Proposition}
\newtheorem{kohilf}[dehilf]{Corollary}
\newtheorem{thhilf}[dehilf]{Theorem}

\newenvironment{De}{\begin{dehilf}
    \hspace{-1.0ex}{\bf :}\begin{rm}}{\end{rm}\end{dehilf}}
\newenvironment{Le}{\begin{lehilf}
    \hspace{-1.0ex}{\bf :}}{\end{lehilf}}
\newenvironment{Sa}{\begin{sahilf}
    \hspace{-1.0ex}{\bf :}}{\end{sahilf}}
\newenvironment{Ko}{\begin{kohilf}
    \hspace{-1.0ex}{\bf :}}{\end{kohilf}}
\newenvironment{Th}{\begin{thhilf}
    \hspace{-1.0ex}{\bf :}}{\end{thhilf}}
\def\qedbox{\mbox{\large $\Box$}}

\newcommand{\deflabel}[1]{\bf1\hfill}%
{\begin{list}{}%
{\settowidth{\labelwidth}{\bf1}%
\setlength{\leftmargin}{\labelwidth}%
\addtolength{\leftmargin}{\labelsep}%
}}%
{\end{list}}%

\newenvironment{Bw}{\begin{list}{}%
       {\setlength{\leftmargin}{0mm}}%
       \item {\bf Proof:}
             }{\rule{0mm}{0mm}\nolinebreak\hfill\qedbox
        \end{list}}

\begin{document}

\title{Codistances of 3-spherical buildings}

\author{Alice Devillers\thanks{Most of the work for this paper was done while the first author was Collaborateur Scientifique of the Fonds National de la
Recherche Scientifique (Belgium); e-mail: \texttt{alice.devillers@uwa.edu.au}}\\
\and Bernhard M\"uhlherr\thanks{e-mail:  \texttt{Bernhard.Muehlherr@math.uni-giessen.de}}\\
\and Hendrik Van Maldeghem \thanks{e-mail:  \texttt{hvm@cage.UGent.be}}}
\date{}\maketitle
\vspace*{-1cm}
%\begin{center}Universit\'e Libre de Bruxelles\\
%D\'epartement de Math\'ematiques - C.P.216\\
%Boulevard du Triomphe\\
%B--1050 Brussels, Belgium
%\end{center}

\begin{abstract} We show that a 3-spherical building in which each rank 2 residue is connected far away from a
chamber, and each rank 3 residue is simply 2-connected far away from a chamber, admits a twinning (i.e., is one
half of a twin building) as soon as it admits a codistance, i.e., a twinning with a single chamber.
\end{abstract}

\section{Introduction}

Twin buildings have been introduced by M.~A~Ronan and J.~Tits in the late 1980's. Their definition is motivated
by the theory of Kac-Moody groups over fields. Kac-Moody groups are infinite-dimensional generalizations of Chevalley groups
and the buildings associated with the latter are spherical. Spherical buildings have been classified by J.~Tits
in \cite{Ti74}. This classification relies heavily on the fact that there is an opposition relation on the set
of chambers of a spherical building. The idea in the definition of a twin building is to extend the notion of
an opposition to non-spherical buildings: instead of taking one building, one starts with two buildings ${\cal
B}_+,{\cal B}_-$ of the same type and defines an opposition relation between the chambers of the two buildings
in question. Technically, this is done by requiring a {\sl twinning function} between the chamber sets of the two
buildings that takes its values in the Weyl group $W$. Two chambers $x,y$ of ${\cal B}_+$ and ${\cal B}_-$ are
then defined to be opposite, if their twinning is the identity in $W$.

There are variations of the idea of a twinning. For instance, one can introduce `by restriction' a twinning
between one chamber of  ${\cal B}_+$ and the building ${\cal B}_-$, seen as an application from the set of chambers
of ${\cal B}_-$ to the Weyl group. A function from the set of chambers of a building ${\cal B}$ to its Weyl group and satisfying similar properties to those of this `twinning to a chamber'  will be 
called {\sl a codistance on ${\cal B}$}.  This idea occurs at various
places in the literature (see for instance \cite{Mu98} and \cite{Ro07}). In particular, \cite{Ro07}  deals
with the question to which extent the existence of a codistance of a building ${\cal B}$ restricts its
structure. The main result of the present paper ensures that any 3-spherical building admitting a codistance,
and satisfying some local condition, is in fact one `half' of a twin building. It is known that the local condition in question 
is satisfied  if the diagram is simply laced and if each panel contains at least 4 chambers (see the final section of this paper).

Here is the precise statement of our main result. For the definitions and
notation we refer to Sections~\ref{prel} and~\ref{cod}.

\smallskip
\noindent {\bf Main result:} \emph{Let ${\cal B}_- =({\cal C}_-,\delta_-)$ be a thick building of 
3-spherical
type $(W,S)$. Assume that the following two conditions hold.
\begin{itemize}
\item[{\rm (lco)}] If $R$ is a rank 2 residue of ${\cal B}_-$ containing a chamber $c$, then the chamber system defined by the set of chambers opposite
$c$ inside $R$ is connected. \item[{\rm (lsco)}] If $R$ is a rank 3 residue  of ${\cal B}_-$ containing a chamber $c$, then the chamber system defined by
the set of chambers opposite $c$ inside $R$ is simply 2-connected.
\end{itemize}
If there exists a codistance function $f: {\cal C}_- \rightarrow W$, then there exists a building ${\cal B}_+
=({\cal C}_+,\delta_+)$ and a mapping $\delta_*: ({\cal C}_- \times {\cal C}_+) \cup ({\cal C}_+ \times {\cal
C}_-) \rightarrow W$ such that the following two statements hold.
\begin{itemize}
\item[{\rm a)}] $({\cal B}_-,{\cal B}_+,\delta_*)$ is a twin building. \item[{\rm b)}] There exists a
chamber $c \in {\cal C}_+$ such that $\delta_*(c,x) = f(x)$ for all $x \in {\cal C}_-$.
\end{itemize}
}

We would like to mention that the Conditions (lco) and (lsco) are `almost always'  automatic in 
3-spherical buildings.
We will explain this in the final section of this paper. In view of the discussion there the following 
corollary is a consequence of our main result.

\smallskip
\noindent {\bf Corollary~1:} \emph{Let ${\cal B}_- =({\cal C}_-,\delta_-)$ be a thick, irreducible building of 3-spherical
type $(W,S)$ whose rank is at least 3.
Then the conclusions of the main result hold as soon as one of the following conditions
is satisfied: 
\begin{itemize}
\item[{\rm (1)}] $(W,S)$ is simply laced and all panels contain at least 4 chambers.
\item[{\rm (2)}] Any residue of type $A_2$ corresponds to a Desarguesian projective plane 
and any panel contains at least 17 chambers. 
\end{itemize}
}

\subsection*{Some general remarks on 3-spherical buildings}

The most impressive results in the theory of abstract buildings are the classifications of the irreducible 
spherical buildings of rank at least 3 and the irreducible affine buildings of rank at least 4 
by Tits. In the 1980's it was an open question whether Tits' classification could be extended
to irreducible affine buildings of rank 3. By independent work of Ronan and the third author constructions
for such buildings were given, which showed that such a classification cannot be expected.
Especially, Ronan's construction could be extended in order to show that there is a sort of 
free construction for buildings of type $(W,S)$ for a lot of Coxeter systems (this is the Ronan-Tits construction, see \cite{RT87}). However,
if a Coxeter system  contains spherical subsystems of rank 3, 
the degree of freeness of the Ronan-Tits-construction is considerably reduced. In fact, if all rank 3 subsystems
are spherical, the only known choice of parameters in that construction yields buildings coming
from Kac-Moody groups---hence buildings of algebraic origin. Note that an irreducible affine building
of rank at least 4 is 3-spherical and that there are affine buildings which do not come from Kac-Moody groups.
At present, the only known irreducible 3-spherical buildings of rank at least 4 which are not of affine
type are those coming from groups of Kac-Moody-type. Whether these are all, appears to be an interesting open question
in the theory of abstract buildings. Our main result provides a positive answer to that question
under the additional assumption that the building admits a codistance. We make this more precise
for the \emph{simply laced} case  (i.e. if all entries are
2 or 3).

\subsubsection*{Simply laced 3-spherical buildings}  
Let $\mathcal{B} = (\mathcal{B}_+,\mathcal{B}_-,\delta_*)$ be an irreducible %3-spherical
twin building of rank at least 3 whose diagram is simply laced. Then it is known that $\mathcal{B}$ is Moufang
(see for instance \cite{Abr-Bro}) and therefore each of its spherical residues is Moufang.
If we assume in addition that $\mathcal{B}$ is 3-spherical, then all its $A_2$-residues
are (up to duality) isomorphic to the building associated to a projective plane over
a division ring $K$. If there is a $D_4$-subdiagram, then $K$ is commutative and those buildings have been classified in \cite{Mu99a}; in particular, they are of
`algebraic origin' in the sense that they can be constructed as `$k$-forms' of certain Kac-Moody groups.
If there is no $D_4$-subdiagram, then $\mathcal{B}$ is of type 
$A_n$ or $\tilde{A}_n$ for some $n \geq 3$. Those buildings
% are spherical of rank at least 3 or affine of rank at least 4 respectively. Therefore, they
are also known by \cite{Ti74} and \cite{Ti84} and of algebraic origin. Putting together
all this information, we get the following corollary of our main result.

\smallskip
\noindent {\bf Corollary~2:}  \emph{Let $\mathcal{B}_-$ be an irreducible, 3-spherical and simply laced
building of rank at least 3 in which each panel contains at least 4 chambers. If $\mathcal{B}_-$ admits
a codistance, then it is known and in particular of algebraic origin. 
}

%\subsubsection*{Condition (H)}
%
%As mentioned above, the question to which extent  the existence of a codistance on 
%a building restricts its structure has been also considered 
%in  \cite{Ro07}. The results obtained there are proved under the additional
%assumption that the building in question satisfies a certain additional condition. In  loc.\,cit.\, 
%it is called Condition (H).  Condition (H) does not play any role
%in our paper. Nevertheless, it is perhaps worthwhile  to mention
%that our main result implies that the buildings 
%covered in its statement do satisfy Condition (H).

\subsection*{Content}

The paper is organized as follows. In Section~\ref{prel}, we collect the definitions, known results and
preliminaries that we need. In Section~\ref{cod}, we prove some basic properties of a codistance; most
properties are known to be valid for a twinning, but we need to reprove them here for a codistance. In
Section~\ref{simpleconnectivity}, we show that, under the assumptions of our main result, the complex of
chambers with codistance the identity is simply 2-connected (for \emph{any} codistance!). In
Section \ref{sec:parpanels} we study parallel panels and in
Section~\ref{bijections}, we construct bijections between panels that are contained in a chamber of codistance
the identity. These bijections will then be used in Section~\ref{adjacentco} to define codistances adjacent to
a given codistance. Finally, in Section~\ref{constructiontwinning}, we prove that all the codistances thus
obtained constitute the second half of a twinning, the first half of which is the original building.

\bigskip
\bigskip
\noindent
{\bf Acknowledgement:} We thank the referee for his detailed comments
and valuable suggestions which improved the presentation of our work considerably.

\section{Preliminaries}\label{prel}
In this section, we recall basic definitions and results.

\subsection*{Chamber systems}

Let $I$ be a set. A {\it chamber system} over $I$ is a pair
${\cal C}=(C, (\sim_i)_{i \in I})$ where $C$ is a set
whose elements are called {\it chambers} and where
$\sim_i$ is an equivalence relation on the set of chambers
for each $i \in I$, such that if $c \sim_i d$ and $c \sim_j d$ then either $i=j$ or $c=d$.

We refer to \cite{Abr-Bro,DM07} for the definitions of {\it $i$-adjacent chambers}, {\it galleries}, {\it
$J$-galleries} {\it $J$-residues}, {\it $i$-panels}. The $J$-residue containing the chamber $c$ is denoted by $R_J(c)$.

Two galleries $G =(c_0,\ldots,c_k)$ and
$H = (c_0',\ldots,c_{k'}') $ with $c_0=c'_0$ and $c_k=c'_{k'}$ are said to be
{\it elementary $2$-homotopic} if
there exist two galleries $X,Y$
and two $J$-galleries $G_0,H_0$ for some $J \subset I$
of cardinality at most $2$ such that $G = XG_0Y$, $H=XH_0Y$.
Two galleries $G,H$ are said to be {\it $2$-homotopic} if
there exists a finite sequence $G_0,G_1,\ldots,G_l$ of  galleries
such that $G_0 = G,G_l = H$ and such that
$G_{\mu-1}$ is elementary $2$-homotopic to $G_{\mu}$ for all
$1 \leq \mu \leq l$.
The chamber system ${\cal C}$ is called
{\it simply $2$-connected} if it is connected and if each closed
gallery is $2$-homotopic to a trivial gallery.

\subsection*{Coxeter systems}

A {\it Coxeter system} is a pair $(W,S)$ consisting of a group $W$ and
a set $S \subset W$ such that $\langle S \rangle = W$,
$s^2 = 1_W \neq s$ for all $s \in S$ and such that the set $S$
and the relations $( (st)^{o(st)})_{s,t \in S}$ constitute
a presentation of $W$, where $o(g)$ denotes the order of $g$.

Let $(W,S)$ be a Coxeter system. The matrix $M(S) := (o(st))_{s,t \in S}$ is called the {\it type} or the {\it
diagram} of $(W,S)$. For an element $w \in W$ we put $l(w):= \min \{ k \in {\bf N} \mid w = s_1s_2 \ldots s_k
\mbox{ where }s_i \in S \mbox{ for } 1 \leq i \leq k \}$. The number $l(w)$ is called the {\it length} of $w$.
For a subset $J$ of $S$ we put $W_J := \langle J \rangle$ and we call it {\it spherical} if $W_J$ is finite.

The following proposition  collects several basic facts on Coxeter groups. 
% which can be found in the usual standard references \cite{Bo68} or \cite{Hu90}. 
These facts will be used without reference 
throughout the paper.

\begin{Sa} \label{sa41}
Let $(W,S)$ be a Coxeter system.
\begin{itemize}
\item[a)] For $w \in W, s \in S$ we have
 $\{l(ws), l(sw) \} \subset \{ l(w)-1, l(w) + 1 \}$.
\item[b)] For $w \in W, s,t \in S$ with $l(sw) = l(w)+1 = l(wt)$ we have
$l(swt) = l(w)+2$ or $swt=w$.
\item[c)] For $J \subset S$ the pair $(W_J,J)$ is a Coxeter system
and if $l_J: W_J \rightarrow {\bf N}$ is its length function, then
$l_J = l \mid_{W_J}$.
\item[d)] Let $w \in W$ and $J \subset S$. Then there exists a unique
element $w_J \in wW_J$ such that $l(w_Jt) = l(w_J) + 1$
for all $t \in J$. Moreover, we have $l(x) = l(w_J) + l_J(w_J^{-1}x)$ for all
$x \in wW_J$.
\item[e)] If $ J \subset S$ is spherical, then there is a unique element $r_J \in W_J$ such that
$l(r_Jw) + l(w) = l(r_J)$ for all $w \in W_J$;
the element $r_J$ is a non-trivial involution if $J \neq \emptyset$.
Moreover, we have $r_JJr_J=J$.
%\item[f)] For each $w \in W$ the set $J^-(w) := \{ t \in S \mid l(wt) = l(w)-1 \}$ is a spherical
%subset of $S$ and $l(wr_{J^-(w)}) = l(w)-l(r_{J^-(w)})$.
\item[f)] Let $w \in W$ and let $J \subset S$ be spherical. Then there exists
a unique element $w^J \in wW_J$ such that $l(w^Jt) = l(w^J)-1$ for all $t \in J$ and
we have $w^J = w_J r_J$.
Moreover we have $l(x) = l(w^J)-l_J((w^J)^{-1}x)$
for all $x \in wW_J$; in particular,
$l(w_J)+l(r_J) = l(w^J)$.
%\item[g)] Let $w \in W$ and $K,J \subset S$, then there exists a unique
%shortest element $\omega$ in $W_KwW_J$. Moreover, if
%$w' \in W_KwW_J$ then $w' = \omega$ if and only if
%$l(w't) = l(w')+1 = l(uw')$ for all $t \in J$ and all $u \in K$.
\end{itemize}
\end{Sa}

\begin{Bw}
Parts a) and b) follow from the Deletion and the Folding Condition (see Sections 2.1 and 2.3 in \cite{Abr-Bro})
and Part c) is an immediate consequence of the Deletion Condition. Part d) is a reformulation of \cite[Proposition 1.10]{Hu90} and
Part e) follows from  \cite[Proposition 5.7]{We03}. Finally, Part f) is a consequence of Parts d) and e).
\end{Bw}

Let $(W,S)$ be a spherical Coxeter system and let $r := r_S$ be the longest element in $W$. Then $rSr = S$ and hence
conjugation by $r$ induces an involutory permutation $\op$ of $S$; we say that $J \subset S$ is opposite $K \subset S$
if $\op(J) = K$. More generally, if $(W,S)$ is an arbitrary Coxeter system and if $J \subset S$ is a spherical subset,
then we say that two sets $K,L \subset J$ are opposite with respect to $J$ if $r_JKr_J = L$.

\subsection*{Buildings}

Let $(W,S)$ be a Coxeter system.
A {\it building} of type $(W,S)$ is a pair ${\cal B}
=(C,\delta)$ where $C$ is a set
and where
$\delta : C \times  C \rightarrow W$ is a {\it distance
function}
satisfying the following axioms, where $x,y \in  C$ and
$w = \delta(x,y)$:

\begin{itemize}
\item[(Bu 1)] $w = 1$ if and only if $x = y$;
\item[(Bu 2)] if $z \in  C$ is such that $\delta(y,z) = s \in S$,
then $\delta(x,z) = w$ or $ws$, and if, furthermore, $l(ws) = l(w) + 1$,
then $\delta(x,z) = ws$;
\item[(Bu 3)] if $s \in S$, there exists $z \in  C$ such that
$\delta(y,z) = s$ and $\delta(x,z) = ws$.
\end{itemize}

For a building ${\cal B} = ( C,\delta)$ we define the chamber
system ${\bf C}({\cal B}) = ( C,(\sim_s)_{s \in S})$ where
two chambers $c,d \in  C$ are defined to be \textit{$s$-adjacent}
if $\delta(c,d) \in \langle s \rangle$. The \textit{rank} of a building ${\cal B}$ of type  $(W,S)$ is $|S|$.

In this paper all buildings are assumed to be of finite rank and
{\it thick} (which means that for any $s \in S$ and any chamber $c \in C$
there are at least three chambers being $s$-adjacent to $c$). 

For any two chambers $x$ and $y$ we set $l(x,y)=l(\delta(x,y))$.
We say that a gallery $x_0,x_1,\ldots,x_n$ is {\it minimal} if $n=l(x_0,x_n)$.

In the following proposition we collect several basic facts about buildings.
%We refer to \cite{Ro89} and \cite{We03} for the details.

\begin{Sa} \label{sa42}
Let $(W,S)$ be a Coxeter system and let ${\cal B}
=(C,\delta)$ be a building of type $(W,S)$.
\begin{itemize}
\item[a)] The chamber system  ${\bf C}({\cal B}) = ( C,(\sim_s)_{s \in s})$ uniquely determines ${\cal B}$; in
other words, the $s$-adjacency relations on $C$ determine the distance function $\delta$. \item[b)] For $c \in
C$ and $J \subset S$ we have
 $R_J(c) = \{ x \in C \mid \delta(c,x) \in W_J \}$.
\item[c)] If $d: C \times C \rightarrow {\bf N}$ is
the numerical distance between two chambers
in $( C,(\sim_s)_{s \in s})$, then $d = l$.
\item[d)] Let $c \in C$ and let $R \subset C$ be a $J$-residue for some
$J \subset S$. Then there exists a unique chamber $x \in R$ such that
$\delta(c,x) = (\delta(c,x))_J$. Moreover, for all $y \in R$ one has
$\delta(c,y) = \delta(c,x) \delta(x,y)$ and in particular,
$l(c,y) = l(c,x) + l(x,y)$.
\end{itemize}
\end{Sa}

\begin{Bw}
Parts a), b) and c) follow from the fact that the distance function can be characterized
in terms of types of galleries (see   \cite[Definition 7.1]{We03} or \cite{Ti81}).
For Part d) we refer to  \cite[Proposition 8.24]{We03}.
\end{Bw}

Given $c \in
 C$ and a $J$-residue $R$ of ${\cal B}$ as in Assertion
d) of the previous proposition,
 then the chamber $x$ in that statement
is called the {\it projection of $c$ onto $R$} and it is denoted by $\proj_R c$.

Given two residues $R$ and $R'$, we define $\proj_{R}R'$ by the set $\{\proj_R c|c\in R'\}$.

Two residues $R_1$ and $R_2$ of a building are called {\it parallel} if $\proj_{R_1}:R_2 \rightarrow R_1$ and
$\proj_{R_2}:R_1 \rightarrow R_2$ are adjacency-preserving bijections inverse to each other.

\begin{Sa}\label{Sa46}
Let $R,Q$ be two residues of a  building.
%Put $\proj_R Q := \{ \proj_R x \mid x \in Q \}$.
Then the following holds:
\begin{itemize}
\item[a)] $\proj_R Q$ is a residue contained in $R$.
\item[b)]  The residues $ R' := \proj_R Q$ and $Q' := \proj_Q R$ are parallel.
\end{itemize}
\end{Sa}
\begin{Bw}
Part a) follows from the first statement of  \cite[Proposition 3]{DS}
while Part b) follows from Part a of the main Theorem in \cite{DS}.
\end{Bw}

Let $R$ be a spherical $J$-residue of a building of type $(W,S)$. Two chambers $x,y$ of $R$ are {\it opposite}
in $R$ whenever $\delta(x,y)=r_J$. Two residues $R_1$ of type $K_1$ and $R_2$ of type $K_2$ in $R$ are {\it
opposite} in $R$ if $R_1$ contains a chamber opposite to a chamber of $R_2$ and if $K_1=r_JK_2r_J$ (which means
that $K_1$ and $K_2$ are opposite with respect to $J$ as defined earlier).

%What's a reference for the following statements?

\begin{Sa}\label{Sa47}
Let $R$ be a spherical $J$-residue of a building of type $(W,S)$ and let $R_1,R_2$ be two residues which are opposite in
$R$. Then $R_1$ and $R_2$ are parallel.
\end{Sa}

\begin{Bw}
This is a consequence  of Theorem 3.28 of \cite{Ti74}.
\end{Bw}

\begin{Sa}\label{Sa48} Let $R_I,R_J,R_K$ be residues of respective type $I,J,K$ of a building  of type $(W,S)$.
Assume that $R_I\subseteq R_J$. Then we have $\proj_{R_I}R_K=\proj_{R_I}\proj_{R_J}R_K$.
\end{Sa}

\begin{Bw}
This is a consequence of Proposition 2 in \cite{DS}
\end{Bw}

\subsection*{Twin buildings}

Let ${\cal B}_+ = ( C_+,\delta_+),{\cal B}_- = ( C_-,\delta_-)$ be two buildings of the same type $(W,S)$,
where $(W,S)$ is a Coxeter system. A {\it twinning} between ${\cal B}_+$ and ${\cal B}_-$ is a mapping
$\delta_*: ( C_+ \times  C_-) \cup (C_- \times C_+) \rightarrow W$ satisfying the following axioms, where
$\epsilon \in \{ +,- \}, x \in C_{\epsilon}, y \in  C_{-\epsilon}$ and $w = \delta_*(x,y)$:

\begin{itemize}
\item[(Tw 1)] $\delta_*(y,x) = w^{-1}$;
\item[(Tw 2)] if $z \in  C_{-\epsilon}$ is such that
$\delta_{-\epsilon}(y,z)=s \in S$ and $l(ws) = l(w) - 1$,
then $\delta_*(x,z) = ws$;
\item[(Tw 3)] if $s \in S$, there exists $z \in  C_{-\epsilon}$
such that $\delta_{-\epsilon}(y,z) = s$ and $\delta_*(x,z) = ws$.
\end{itemize}

A {\it twin building of type $(W,S)$} is a triple
$({\cal B}_+,{\cal B}_-,\delta_*)$ where ${\cal B}_+,{\cal B}_-$
are buildings of type $(W,S)$ and where $\delta_*$ is a twinning
between ${\cal B}_+$ and ${\cal B}_-$.

Let ${\cal B} = ({\cal B}_+,{\cal B}_-,\delta_*)$ be a twin building. Then $x \in {\cal C}_+$ and
$y \in {\cal C}_-$ are called \textit{opposite} if $\delta_*(x,y) = 1_W$. For each chamber $c$ in one
of the two buildings, $c^{\op}$ denotes the set of chambers in the other building that are 
opposite $c$.

%We recall that in this paper  ${\cal B}_+$ and ${\cal B}_-$ are thick and of finite rank.

Here is a lemma the proof of which is left to the reader (it follows directly from the definition above and an
easy induction on the length of $\delta_+(x,y)$).

\begin{Le}\label{+*} Let $(({\cal C}_+,\delta_+),({\cal C}_-,\delta_-),\delta_*)$ be a twin building of type $(W,S)$.
Let $x,y \in {\cal C}_+$ and $z \in {\cal C}_-$ be such that $\delta_+(x,y) = \delta_*(x,z)$. Then $y$ and $z$
are opposite. In particular, if $x^{\op} = y^{\op}$, then $x = y$.
\end{Le}

\section{Codistances}\label{cod}

In this section, we take  $(W,S)$ a Coxeter system and ${\cal B}
=({\cal C},\delta)$ a building of type  $(W,S)$.

\begin{De}
A \emph{codistance} on  ${\cal B}$ is a function $f:{\cal C}\rightarrow W$ such that, for all $s \in S$ and $P$
an $s$-panel of ${\cal C}$, there exists $w\in W$ with $f(x)\in \{w,ws\}$ for all $x\in P$ and $P$ contains a
unique chamber with $f$-value the longer word of the two. If the latter is satisfied for a fixed panel $P$, then we say that \emph{$P$ satisfies the codistance condition for $f$}. 
\end{De}

As an example, if ${\cal B}$ is half of a twin building and $x$ is a chamber in the other half, the twinning to
$x$ is a codistance on ${\cal B}$.

For the rest of this section $f:{\cal C} \rightarrow W$ is a codistance of ${\cal B}$.

\begin{Le}\label{lem:imfR}
Let $R$ be a $J$-residue of ${\cal B}$ and $x$ be a chamber of $R$. Then the image of $f$ restricted to $R$ is
$f(x)W_J$.
\end{Le}
\begin{Bw}
By the definition of $f$, the image is contained in $f(x)W_J$. Let $w$ be a word of $W_J$
written as a reduced word as $s_1s_2\ldots s_k$. Using the fact that for all $s\in J$ and all chambers $y\in R$,
there exists at least one chamber $s$-adjacent to $y$ with $f$-value $f(y)s$, it follows by induction
on $k$ that there exists a chamber in $R$ with $f$-value $f(x)w$.
\end{Bw}

\begin{Sa}\label{projf}
Let $R$ be a spherical $J$-residue  of ${\cal B}$. Then there exists a unique chamber $c$ in $R$ such that $l(f(c))>l(f(y))$
for all $y\in R\setminus\{c\}$. This unique chamber will be denoted by $\proj_Rf$. Moreover for all $y\in R$, we
have $f(y)=f(c)\delta(c,y)$.
\end{Sa}
\begin{Bw} Let $y$ be a chamber in $R$ and $w:=f(y)$. By Lemma \ref{lem:imfR}, $f$ takes on $R$ its values in $wW_J$.
Since $R$ is spherical, 
$wW_J$ contains a unique longest word $w^J$ by Part f) of Proposition \ref{sa41}. Moreover $l(x)=l(w^J)-l((w^J)^{-1}x)$ for all $x\in wW_J$. By
Lemma \ref{lem:imfR}, there exists a chamber $c\in R$ with $f(c)=w^J$.

Let $y$ be a chamber in $R$. The distance $\delta(c,y)$ is in $W_J$ and so can be written as a reduced word as
$t_1t_2\ldots t_k$. Therefore there is a gallery  $c= y_0 \sim_{t_1} y_1 \sim_{t_2} \ldots \sim_{t_k} y_k = y$.
Using the fact that $l(w^Jt_1\ldots t_i)=l(w^J)-i$, it follows by induction that
$f(y)=w^J\delta(c,y)$. Therefore $l(f(y))=l(w^J)-l(\delta(c,y))=l(f(c))-l(c,y)\leq l(f(c))$ with equality only
if $y=c$.
\end{Bw}

\begin{Sa} \label{AfR}
Let $R$ be a $J$-residue of ${\cal B}$. Put $l_f(R) := \min \{l(f(x)) \mid x \in R \}$
and $A_f(R) := \{ x \in R \mid l(f(x)) = l_f(R) \}$.
\begin{itemize}
\item[a)] Let $x \in R$. Then $x \in A_f(R)$ if and only if
$f(x)$ is the unique shortest word of $f(x)W_J$. Moreover, if $x,y \in A_f(R)$,
then $f(x) = f(y)$.
\item[b)] Let $y \in R$. Then there exists $x \in A_f(R)$,
such that $f(y) = f(x) \delta(x,y)$.
\item[c)] If $J$ is spherical, then $A_f(R)$ is the set of all chambers opposite $\proj_Rf$ in $R$.
%\item[d)] If $J$ is a spherical subset of $S$, then any two $J$-residues
%of ${\cal B}$ are isomorphic.
\end{itemize}
\end{Sa}

\begin{Bw}
Let $y \in R$ and put $w:=f(y)$.

%Suppose first that there is $t \in J$ with $l(wt) = l(w)-1$.
%Choose a chamber $y' \in R$ with $y \sim_t y' \neq y$.
%It follows from the definition of $f$ that $f(y') = wt$.
%Hence $y \not \in A_f(R)$ in this case.

By Lemma \ref{lem:imfR}, $\{ f(x) \mid x \in R \} = wW_J$. By Part d) of Proposition \ref{sa41}
there exists a unique shortest element $w_J\in wW_J$.
It follows that
$A_f(R) = \{ x \in R \mid f(x) = w_J \}$. This proves
Part a) of the proposition.

Now let $t_1 t_2 \ldots t_k$ be a reduced representation of $w_J^{-1}w$ and let $x = y_0 \sim_{t_1} y_1
\sim_{t_2} \ldots \sim_{t_k} y_k = y$ be a reduced gallery ending in $y$. Using the fact that
$l(wt_kt_{k-1}\ldots t_{i+1})=l(w_Jt_1t_2\ldots t_i)=l(w_J)+i$, it follows by induction on $k$ that
$x \in A_f(R)$. By construction $\delta(x,y) =t_1 t_2 \ldots t_k=w_J^{-1}w=f(x)^{-1}f(y)$. This finishes Part
b).

Let $J$ be spherical. 
%Then $x \in A_f(R)$ if and only if $f(x)=w_J$. 
Let $c=\proj_Rf$ so that $f(c)=w^J$ as in
Proposition \ref{projf}. We have seen that $f(x)=w^J\delta(c,x)$ for all $x\in R$. Since $w^J=w_Jr_J$, where $r_J$ is the unique
longest word of $W_J$, we can conclude that $x \in A_f(R)$ if and only if $\delta(c,x)=r_J$, that is, if and
only if $x$ is opposite $c$ in $R$.

%Part c) follows from (4.1) in \cite{Ro00} and Part d) is an easy
%consequence of (4.3) in loc. cit..
\end{Bw}

\begin{De} We denote by $f^{\op}$ the set of chambers of ${\cal C}$ with $f$-value $1_W$.
 \end{De}
 
Let $R \subset {\cal C}$ be a $J$-residue. By Part a) of the previous proposition
we have $f(x) = f(y)$ for all $x,y \in A_f(R)$. We denote this common value by $R^f$. 
Note that $A_f(R) = \{ x \in R \mid f(x)_J = f(x) \}$ by Part a) of the previous proposition.

\begin{Le}\label{galltofop}
Let $c$ be a chamber of ${\cal C}$. Then a shortest gallery from $c$ to a chamber in $f^{\op}$ has length
$l(f(c))$.
\end{Le}
\begin{Bw}
It is obvious from the definition of the codistance $f$ that no chamber at distance strictly less than  $l(f(c))$ from
$c$ can be in $f^{\op}$. Now by Part b) of Proposition \ref{AfR} with $J=S$, there exists $x\in A_f({\cal C})=f^{\op}$
such that $f(c)=\delta(x,c)$. Hence a minimal gallery from $c$ to $x$ will have length $l(f(c))$.
\end{Bw}

For $c\in {\cal C}$, we define $f^{\op}_c=\{x\in f^{\op}|\delta(x,c)=f(c)\}$ (that is the set of chambers 
of $f^{\op}$ closest to $c$), which is non-empty, by Lemma \ref{galltofop}.
\begin{Le}\label{fopc} The following statements are equivalent:
\begin{itemize}
\item[a)] the  chamber $x$ is in $f^{\op}_c$,
\item[b)] for any minimal gallery $x=x_0,x_1\ldots,x_n=c$ we have $l(f(x_i))=i$ for all $0\leq i\leq n$,
\item[c)] there exists a minimal gallery $x=x_0,x_1\ldots,x_n=c$ with $l(f(x_i))=i$ for all $0\leq i\leq n$.
\end{itemize}
\end{Le}
\begin{Bw}
Assume $x\in f^{\op}_c$. Let  $x=x_0,x_1\ldots,x_n=c$ be any minimal gallery from $x$ to $c$. Since
$\delta(x,c)=f(c)$, $n=l(f(c))$. By the axioms of codistance, the $f$-values of two adjacent chambers are
either equal or have length difference one, hence we must have  $l(f(x_i))=l(f(x_{i-1}))+1$ for all $1\leq
i\leq n$, which implies b).

Obviously b) implies c).

Assume that there exists a minimal gallery $x=x_0,x_1\ldots,x_n=c$ with $l(f(x_i))=i$ for all $0\leq i\leq n$.
Then $l(f(x))=0$, so $f(x)=1_W$. Assume that $f(x_i)=\delta(x,x_i)$, then  $f(x_{i+1})=\delta(x,x_{i+1})$.
Indeed $x_i\sim_{s_i}x_{i+1}$ for some $s_i\in S$ and so $f(x_{i+1})=f(x_i)$ or $f(x_i)s_i$. Since
$l(f(x_{i+1}))\neq l(f(x_i))$, we are in the second case and  $f(x_{i+1})=\delta(x,x_i)s_i=\delta(x,x_{i+1})$.
This proves by induction that $f(x_i)=\delta(x,x_i)$ for all $0\leq i\leq n$, and so $f(c)=\delta(x,c)$, which
yields a).
\end{Bw}

\begin{Le}\label{uniquec}
Let $x\in{\cal C} $ and $w\in W$ such that $l(f(x)w)=l(f(x))+l(w)$. Then there exists a unique chamber $c$ of
${\cal C}$ with $f(x)^{-1}f(c)=w=\delta(x,c)$.
\end{Le}
\begin{Bw}
Let $s_1s_2\ldots s_k$ be a reduced word for $w$. Since  $l(f(x)w)=l(f(x))+l(w)$, we have $l(f(x)s_1s_2\ldots
s_i)=l(f(x))+i$. Consider the $s_1$-panel on $x$, it follows from the axioms of
codistance that this panel contains a unique chamber with $f$-value $f(x)s_1$, namely the projection of $f$ on
it. Continuing by induction on $k$, we can build a unique gallery $x=x_0\sim_{s_1}x_1\sim{s_2}\ldots
\sim_{s_k}x_k=c$ such that $f(x_i)=f(x)s_1s_2\ldots s_i$ for all $i$. Hence $w=\delta(x,c)$ and $f(c)=f(x)w$,
and so $c$ exists.

Assume there exists another chamber $c'$ with $f(x)^{-1}f(c')=w=\delta(x,c')$. Hence, on the one hand, there exists a minimal
gallery $x=x'_0,x'_1\ldots,x'_k=c'$ with $l(f(x_i))=f(x)+i$ for all $0\leq i\leq k$, of type
$t_1,t_2,\ldots,t_k$ where $t_1t_2\ldots t_k=w$. On the other hand, since $\delta(x,c)=w$, there is a minimal
gallery of type  $t_1,t_2,\ldots,t_k$ from $x$ to $c$. Because the length of the $f$-value has to increase at each step, we see by
induction that this gallery coincides with   $x=x'_0,x'_1\ldots,x'_n=c'$, and so $c=c'$.
\end{Bw}

\begin{Le}\label{projAfR}
Let $R$ be a $J$-residue of ${\cal B}$ and $c$ a chamber of $R$. If $x\in f_c^{\op}$ then $\proj_Rx\in A_f(R)$,
$l(x,\proj_Rx)=l_f(R)$ and $\delta(x,\proj_Rx)=R^f$.
\end{Le}
\begin{Bw}
Let $w=f(c)=\delta(x,c)$.
%There exists a minimal gallery from $x$ to $c$ going through $\proj_Rx$: $x=x_0,x_1,\ldots,x_k=\proj_Rx,x_{k+1},\ldots,x_n=c$. By Lemma \ref{fopc}, we have  $l(f(x_i))=i$ for all $0\leq i\leq n$. Therefore, again by  Lemma \ref{fopc}, $x\in f^{\op}_{x_i}$, and so $\delta(x,x_i)=f(x_i)$  for all $0\leq i\leq n$.
%Suppose $x_k\not\in A_f(R)$. By Proposition \ref{AfR}, there exists $y\in  A_f(R)$ with $f(x_k)=f(y)\delta(y,x_k)$. Moreover $l(f(x_k))=k=l_f(R)+l(\delta(y,x_k))$.
%%%
% Then $c$ is the unique chamber with $f(c)=w=\delta(x,c)$.
We have $l(w)=l(w_J)+l(w_J^{-1}w)$. Hence, if  $s_1s_2\ldots s_k$ is a reduced word for $w_J$ and
$s_{k+1}s_{k+2}\ldots s_n$ is a reduced word for $w_J^{-1}w\in W_J$, then $s_1s_2\ldots s_n$ is a reduced word
for $w$. Consider the gallery $x=x_0\sim_{s_1}x_1\sim_{s_2}\ldots\sim_{s_n}x_n$ and such that $f(x_i)=s_1s_2\ldots s_i$.
In particular $f(x_n)=w$ and $f(x_k)=w_J$. By construction, since $s_1s_2\ldots s_i$ is a reduced word, we also have $\delta(x,x_i)=s_1s_2\ldots s_i$ and in particular $\delta(x,x_n)=w$. A
chamber satisfying $ f(x_n)=w=\delta(x,x_n)$ is unique by Lemma \ref{uniquec} and therefore $x_n=c$. As
$w_J^{-1}w\in W_J$, $s_i\in W_J$ for $i\geq k+1$, and so $x_i\in R$  for $i\geq k$. Since
$l(x,x_k)=l(\delta(x,x_k))=l(w_J)$, which is the shortest possible length for under the restriction $x_k\in R$, we have $x_k=\proj_Rx$ and
$x_k\in A_f(R)$ by Proposition \ref{AfR} a). This shows in particular that $\delta(x,\proj_Rx)=R^f$ because
$R^f = w_J$, and so $l(x,\proj_Rx)=l_f(R)$.
\end{Bw}

\begin{Le}\label{fopunique} The set $f^{\op}$ uniquely determines $f$.
\end{Le}
\begin{Bw}
Assume there exists a codistance $f'\neq f$ on  ${\cal B}$ with $f'^{\op}=f^{\op}$. Then consider $c$ at
minimal distance from  $f^{\op}$ under the condition that $f'(c)\neq f(c)$. Of course, $c$ is not in $f^{\op}$.
Let $c=c_0,c_1,\ldots,c_m$ be a shortest gallery from $c$ to $f^{\op}$. This minimal gallery has length
$l(f(c))$ by Lemma \ref{galltofop}. It is also a shortest gallery to $f'^{\op}$, and so has length $l(f'(c))$.
Therefore $l(f(c))=l(f'(c))$. Now $c_1$ is closer to  $f^{\op}$ than $c$, and so $f(c_1)=f'(c_1)$. By the
definition of codistance, $f(c)=f(c_1)$ or $f(c_1)t$ (where $t$ is such that  $c_0\sim_t c_1$). This holds also with  $f'$
in place of $f$.
Since $l(f(c))=l(f'(c))$, it implies that $f(c)=f'(c)$. This contradiction proves that $f=f'$.
\end{Bw}

\section{Simple connectivity of $f^{\op}$}\label{simpleconnectivity}

In this section we will apply a result proved in \cite{DM07} using filtrations.

Let $I$ be a set and let ${\cal C} = (C,(\sim_i)_{i \in I})$ be
a chamber system over $I$.
In the following we denote the set of non-negative integers by
${\bf N}$ and the set of positive integers by ${\bf N}_0$.

A \textit{filtration} of ${\cal C}$ is a family
${\cal F} = (C_n)_{n \in {\bf N}}$
of subsets of $C$ such that the following holds.
\begin{itemize}
\item[(F1)] $C_n \subset C_{n+1}$ for all $n \in {\bf N}$,
\item[(F2)] $\bigcup_{n \in {\bf N}} C_n = C$,
\item[(F3)] for each $n >0$ if $C_{n-1}\neq \emptyset$ then there exists an index $i \in I$
such that for each chamber $c \in C_n$ there exists a chamber
$c' \in C_{n-1}$ which is $i$-adjacent to $c$.
\end{itemize}

A filtration ${\cal F} = (C_n)_{n \in {\bf N}}$ is called \textit{residual} if for each $\emptyset \neq J
\subset I$ and each $J$- residue $R$ the family $(C_n \cap R)_{n \in {\bf N}}$ is a filtration of the chamber
system ${\cal R} := (R, (\sim_j)_{j \in J})$.

For each $x \in C$ we put $|x| := \min\{ \lambda \in {\bf N}
\mid x \in C_{\lambda} \}$. For a subset $X$ of $C$
we put $|X| := \min\{ |x| \mid x \in X \}$ and $\aff(X) :=
\{ x \in X \mid |x| = |X| \}$. Note that $C_0 = \aff(C)$
if we assume that $C_0 \neq \emptyset$.

We say that a filtration \emph{satisfies Condition} (lco)
if for every rank 2 residue $R$, $\aff(R)$ is a connected subset of the chamber system
${\cal R}$.

We say that a filtration \emph{satisfies Condition} (lsco)
if for every rank 3 residue $R$, $\aff(R)$ is a simply 2-connected subset of the chamber system
${\cal R}$.

\begin{thhilf}[see \cite{DM07}]
    \hspace{-1.0ex}{\bf :}\label{th315}
Suppose that the residual filtration ${\cal F}= (C_n)_{n \in {\bf N}}$ of the chamber system $\cal C$
satisfies {\rm (lco), (lsco)} and that $C_0 \neq \emptyset$. Then the following are equivalent:
\begin{itemize}
\item[a)] ${\cal C}$ is simply 2-connected;
\item[b)]  $(C_n, (\sim_i)_{i \in I})$
is simply 2-connected for all $n \in {\bf N}$.
\end{itemize}
\end{thhilf}

\subsection*{The filtration ${\cal F}_f$}

We choose an injection $w \mapsto |w|$ from $W$ into ${\bf N}$
such that $l(x) < l(y)$ implies $|x| < |y|$ for all $x,y \in W$ and such that $|1_W|=0$.
Such an injection exists because ${\cal B}$ is of finite rank.
Let $f$ be a codistance.
We define $C_n$ by setting $C_n := \{ x \in C \mid |f(x)| \leq n \}$.

The goal of this subsection is to show the following proposition.

\begin{Sa}\label{resfilt}
With the definitions above, the family ${\cal F}_f := (C_n)_{n \in {\bf N}}$
is a residual filtration of the chamber system ${\cal C}$.
\end{Sa}

\begin{Bw}
It is obvious that ${\cal F}_f$ satisfies the axioms (F1) and (F2) and from this it follows
that these axioms also hold `residually'.

Let $R$ be a $J$-residue of ${\cal C}$ with $J \neq \emptyset$ and let $|R| := \min \{ k \mid C_k \cap R \neq
\emptyset \}$. It follows from the definition of ${\cal F}_f$ and by Proposition \ref{AfR} that
$\aff(R)=C_{|R|} \cap R = A_f(R) = \{ x \in R \mid f(x) = f(x)_J \}$.

Let $0<n \in {\bf N}$ be such that $C_{n-1} \cap R \neq \emptyset$. We have to show that there is $t \in J$
with the property that each chamber $x$ in $R \cap C_n$ is $t$-adjacent to a chamber $x' \in R \cap C_{n-1}$.
If $C_n \cap R = C_{n-1} \cap R$ we can choose $t \in J$ arbitrarily and set $x' := x$ for each $x \in R \cap
C_n$. Suppose now that $ C_{n-1} \cap R$ is properly contained in $C_n \cap R$, choose $y \in R\cap C_n
\setminus C_{n-1}$ and put $w := f(y)$. Since $|~\cdot~|$ injects $W$ into  ${\bf N}$, it follows from the
definition of ${\cal F}_f$ that, on the one hand, $f(y') = w$ for all $y' \in C_n \setminus C_{n-1}$. On the other hand, there
exists $x \in A_f(R)$ such that $w = f(y) = f(x) \delta(x,y)$ by Assertion b) of Proposition \ref{AfR}. As
$C_{n-1} \cap R \neq \emptyset$ it follows that $y \not \in A_f(R)$ and hence $\delta(x,y) \in W_J \setminus \{
1_W \}$. Let $t \in J$ be such that $l(\delta(x,y)t) = l(\delta(x,y))-1$. As $f(x) = f(x)_J=w_J$ and
$\delta(x,y) \in W_J$, it follows that $l(wt) = l(w_J\delta(x,y)t)=
l(w_J)+l(\delta(x,y)t)=l(w_J)+l(\delta(x,y))-1=l(w)-1$, by Part d) of Proposition \ref{sa41}. For any chamber $z \in R \cap
C_n$ we choose a chamber $z' \in R$ as follows. If $z \in C_{n-1}$ then we put $z' := z$. If $z \in C_n
\setminus C_{n-1}$ then we know that $f(z)=w$ and we choose $z' \in R$ such that $z \sim_t z' \neq z$. In the
first case, it is obvious that $z'$ is in $R \cap C_{n-1}$; in the second case we have $f(z') = wt$ by the
definition of $f$, as $wt$ is shorter than $w$. It follows that $|wt| < |w| = n$ and therefore $z' \in
C_{n-1}$. As $t \in J$ we have also $z' \in R$.

The case $J=S$ is a special case of the consideration above. This
shows that ${\cal F}_f$ satisfies Axiom (F3). Hence ${\cal F}_f$
is a residual filtration.
\end{Bw}

\begin{Th} \label{th52}
Let  ${\cal B}
=({\cal C},\delta)$ be a building of type  $(W,S)$ and $f$ a codistance on ${\cal B}$.
Suppose that the following conditions are satisfied:
\begin{itemize}
\item[{\rm (3-sph.)}] If $J \subseteq S$ is of cardinality at most 3, then $J$ is spherical. \item[{\rm (lco)}] If $J$ is
of cardinality 2, if $R \subset {\cal C}$ is a $J$-residue and if $x \in R$, then the chamber system  $(\{ y
\in R \mid \delta(x,y) = r_J \}, (\sim_t)_{t \in J})$ is  connected. \item[{\rm (lsco)}]  If $J$ is of cardinality 3,
if $R \subset {\cal C}$ is a $J$-residue and if $x \in R$, then the chamber system $(\{ y \in R \mid
\delta(x,y) = r_J \}, (\sim_t)_{t \in J})$ is  simply 2-connected.
\end{itemize}
Then the chamber system $f^{\op}$ is simply 2-connected.
\end{Th}

\begin{Bw}
Let ${\cal F}_f = (C_n)_{n \in {\bf N}}$ be the residual filtration
of Proposition \ref{resfilt}. Note first that $C_0 =f^{\op}$.

Given a spherical $J$-residue $R$ of ${\cal B}$, then
$\aff(R) = A_f(R)$ as we have proved above. By Assertion c) of Proposition \ref{AfR},
 we have therefore
$\aff(R) = \{ x \in R \mid \delta(\proj_R f, x) = r_J \}$, where $r_J$ is the longest
word of $W_J$.

Now ${\cal F}_f$ satisfies (lco)  and (lsco). As it is well-known that ${\cal C}$
is simply 2-connected (see for instance  Theorem (4.3) in \cite{Ro89}), the claim follows now from Theorem \ref{th315}.
\end{Bw}

\section{Parallel panels in buildings}\label{sec:parpanels}

In this section $(W,S)$ is a Coxeter system and ${\cal B} =({\cal C},\delta)$ is a building of type $(W,S)$.

\begin{De}
For $w \in W$ we put $S^-(w) := \{ s \in S \mid l(ws) = l(w)-1 \}$ and $S^+(w) := \{s \in S \mid l(ws) = l(w)+1 \}$.
\end{De}

\begin{Le} \label{le52n}
Let $w \in W$ and $J \subset S$. Then $w = w_J$ if and only if $J \subset S^+(w)$.
\end{Le}
\begin{Bw} This is a consequence of Part d) in Proposition \ref{sa41}.
\end{Bw}

\begin{Le} \label{le53n} 
Let $w \in W$ and $J \subset S$.
 Then the following are equivalent:
\begin{itemize}
\item[a)] $J \subset S^-(w)$,
\item[b)] $J$ is spherical and $l(w) = l(wr_J) +l(r_J)$,
\item[c)] $J$ is spherical and $w = w_Jr_J$,
\item[d)] $J$ is spherical and $w=w^J$.
\end{itemize}
\end{Le}
\begin{Bw}
This is Lemma 2.8 in \cite{Mu92}.
\end{Bw}

\begin{De}
For $w_1,w_2\in W$, we denote $w_1\prec w_2$ if  $l(w_1^{-1}w_2)=l(w_2)-l(w_1)$.
\end{De}

\begin{Le} \label{le55n}
Let $w \in W$ and $J \subset S$. Then $w_J \prec w$ and if $J$ is spherical,
then $w \prec w^J$.
\end{Le}
\begin{Bw}
This follows from Parts d) and f) of Proposition \ref{sa41}.
\end{Bw}

Let $R  \subset {\cal C}$ be a $J$-residue for some subset $J$ of $S$. We recall that
for each chamber $c \in {\cal C}$ there is a unique chamber $\proj_R c$ in $R$
satisfying $\delta(c,\proj_Rc) = \delta(c,\proj_Rc)_J$  and that we 
have $\delta(c,d) = \delta(c,\proj_Rc) \delta(\proj_Rc,d)$ for any $d\in R$. Hence we have
a mapping $\proj_R:{\cal C} \rightarrow R$.

\begin{Le} \label{DSLemma1}
Let $R$ be a residue of ${\cal B}$ and let $c,d \in {\cal C}$.
Then $l(\proj_Rc,\proj_Rd) \leq l(c,d)$.
\end{Le}
\begin{Bw}
This follows from Lemma 1 in \cite{DS}.
\end{Bw}

Let $R,T \subset {\cal C}$ be two residues. We recall that $\proj_RT$ is a residue contained in $R$
and   that the residues $R$ and $T$ are called parallel if
$\proj_RT = R$ and $\proj_TR=T$. By Part b) of Proposition \ref{Sa46} we know that $\proj_RT$ and $\proj_TR$
are parallel residues.
We denote the restriction of $\proj_R:{\cal C} \rightarrow R$
to the residue $T$ by $\proj_R^T$.  

\begin{Le} \label{le57n}
Let $P$ be a panel and let $R$ be a residue. Then $\proj_RP$
is either a singleton or a panel contained in $R$.
In particular, if $|\proj_RP| \geq 2$, then $\proj_RP$ is a panel
contained in $R$.
\end{Le}
\begin{Bw}
We know that $\proj_RP$ is a residue in which any two chambers have
distance at most 1 by Lemma \ref{DSLemma1}. The claim follows.
\end{Bw}

\begin{Le}\label{condpar} Two panels $P_1$ and $P_2$ are parallel if and only if $|\proj_{P_2}P_1| \geq 2$.
\end{Le}
\begin{Bw}
If $P_1$ and $P_2$ are parallel, then $|\proj_{P_2}P_1| =|P_2| \geq 2$
because $\proj_{P_2}^{P_1}$ is a bijection from $P_1$ onto $P_2$. 

Suppose now that $|\proj_{P_2}P_1| \geq 2$. Then $\proj_{P_2}P_1$ is a panel
contained in $P_2$  by the previous lemma and therefore $P_2 = \proj_{P_2}P_1$.
As $\proj_{P_1}P_2$ is parallel to $\proj_{P_2}P_1 = P_2$ by Assertion b) of Proposition \ref{Sa46},
it follows that $|\proj_{P_1}P_2| = |P_2| \geq 2$. Using the same argument as before
we obtain $P_1 = \proj_{P_1} P_2$. Hence $P_1$ and $P_2$ are parallel.
\end{Bw}

%On the other hand, we know that $\proj_{P_2}P_1$ is a residue contained in $P_2$ by Part a) of
%Proposition \ref{Sa46}. Hence, if this set
%contains at least two chambers it follows $\proj_{P_2}P_1 =P_2$, because $P_2$ is a panel.
%Now $\proj_{P_1}P_2 =\proj_{P_1}  \proj_{P_2}P_1$ is a residue contained in $P_1$, which is 
%parallel with $P_2 =  \proj_{P_2}P_1$ by Part b) of Proposition \ref{Sa46}. It follows in particular,
%that $|\proj_{P_1}P_2| geq 2$ and finally $\proj_{P_1}P_2 = P_1$ the same argument as before.
%Hence $P_1$ and $P_2$ are parallel.
%\end{Bw}

\begin{Le}\label{parallelpanels}
Let  $P_1$ and $P_2$ be two parallel panels of type $s_1$ and $s_2$, respectively. Then  $s_2=w^{-1}s_1w$,
where $w:=\delta(x,\proj_{P_2}x)$ does not depend on the choice of $x$ in $P_1$. 

Conversely, if $x$ and $y$ are chambers with 
$\delta(x,y)=w$, where $w$ satisfies  $s_2=w^{-1}s_1w$ and  $l(s_1w)=l(w)+1$, then the $s_1$-panel on $x$ is parallel to the
$s_2$-panel on $y$.
\end{Le}
\begin{Bw}
We know by Part b) of Proposition \ref{sa41} that for $w\in W$ and $s_1,s_2\in S$ such that $l(s_1w)=l(w)+1=l(ws_2)$, we have
$l(s_1ws_2)=l(w)+2$ or $s_1ws_2=w$. The result follows.
\end{Bw}

%The distance $\delta(x,\proj_{P_2}x)$ between two parallel panels $P_1$, $P_2$ will be denoted  by $\delta(P_1,P_2)$.

\begin{De}
For two parallel panels $P_1$ and $P_2$ we put $\delta(P_1,P_2):=\delta(x,\proj_{P_2}x)$, where $x$
is a chamber in $P_1$; by the previous lemma $\delta(P_1,P_2)$ does not depend on the choice of $x \in P_1$.
 \end{De}
 
\begin{De}
 For $s\in S$, we define $X_s:=\{w\in W|w^{-1}sw\in S \text{ and }l(sw)=l(w)+1\}$. 
\end{De}

\begin{Le}\label{Xs}
For $w\in X_s$ and a given $s$-panel $P$, there exists an $w^{-1}sw$-panel $P'$ parallel to $P$ and with
$\delta(P,P')=w$. Let $J$ be a spherical subset of $S$ containing $s$ and let $r_J$ be the longest word of
$W_J$, then $x_J:=sr_J$ is in $X_s$. Moreover, if $w\in W_J$ is in $X_s$, then $w\prec x_J$.
\end{Le}
\begin{Bw}
The first statement is a corollary of Lemma \ref{parallelpanels}. We have $x_J^{-1}sx_J=r_Jsr_J$, which has
length $l(r_J)-l(sr_J)=1$ and so is an element of $S$, and $l(sx_J)=l(r_J)=l(x_J)+1$, hence the second
statement. Finally $l(w^{-1}x_J)=l(w^{-1}sr_J)=l(r_Jsw)=l(r_J)-l(sw)=l(r_J)-l(w)-1=l(x_J)-l(w)$, hence the
third statement.
\end{Bw}

\begin{Le} \label{le512n}
Let $s \in S$, let $w \in X_s$ and put $t := w^{-1}sw \in S$.
Then $t \in S^+(w)$ and $S^-(wt) = S^-(w) \cup \{ t \}$.
In particular $S^-(w) \cup \{ t \}$ is a spherical subset of $J$.
\end{Le}
\begin{Bw}
Note first that $l(w)+1 = l(sw) = l(wt)$ and therefore $t \in S^+(w)$ and $t \in S^-(wt)$.
Now let  $u \in S$ be unequal to $t$. Suppose first that 
$u \in S^-(w)$. Then
$$l(wtu) = l(swu) \leq l(s)+l(wu) = 1 + l(w)-1 = l(w)=l(wt)-1$$
and so we have $l(wtu) = l(wt)-1$ and $u \in S^-(wt)$. Hence $S^-(wt) \supset S^-(w) \cup \{ t \}$.

Suppose now that $u\in S^-(wt)$ but $u\notin S^-(w) \cup \{ t \}$. 
Then $l(wu) = l(w)+1$ and $u\neq t$. By Part b) of Proposition \ref{sa41}, $swu=w$ or $l(swu)=l(w)+2$. 
In the first case, $swu=wtu=w$ implies $t=u$ a contradiction, so, on the one hand,  $l(swu)=l(w)+2$. 
On the other hand $l(wtu) = l(wt)-1=l(w)$, so we also get a contradiction. Therefore  $S^-(wt) \subset S^-(w) \cup \{ t \}$.
This proves the first assertion.
The second assertion is now a consequence of Lemma \ref{le53n}.
\end{Bw}

\begin{De}
Let $\Gamma$ be the graph whose vertices are the panels of ${\cal B}$ with panels adjacent if there exists a
rank 2 residue in which the two panels are opposite. For two adjacent panels $P,Q$, there exists a unique rank
2 residue containing $P$ and $Q$, that will be denoted by $R(P,Q)$. A path $\Pi= (P_0,P_1,\ldots P_k)$  
 in $\Gamma$ is called {\it compatible} if $\proj_{R(P_{i-1},P_i)}P_0=P_{i-1}$ for all $1\leq i\leq
k$. The number $k$ is the {\it length} of that path $\Pi$. The sequence $(J_1,\ldots,J_k)$ where
$J_i$ is the type of $R(P_{i-1},P_i)$ will be called the {\it type} of $\Pi$.
\end{De}

\begin{Le}\label{projpanels}
Let $P,Q$ be two parallel panels of ${\cal B}$ and let $R$ be a residue containing $Q$. Then $\proj_RP$ is a panel
parallel to both $P$ and $Q$. Moreover, if $P=P_0,P_1,\ldots P_k=\proj_RP$ and $\proj_RP=T_0,T_1,\ldots T_l=Q$
are compatible paths in $\Gamma$, the second one contained in $R$, then $P=P_0,P_1,\ldots P_k=T_0,T_1,\ldots
T_l=Q$ is a compatible path in $\Gamma$.
\end{Le}
\begin{Bw}
The projection of a residue on a residue is a residue, so  $P':=\proj_RP$ is either a chamber or a panel. Since
$\proj_Q=\proj_Q\proj_{R}$ by Proposition \ref{Sa48}, we have $Q=\proj_QP=\proj_QP'$, and so $P'$ cannot be reduced to a chamber and is
parallel to $Q$ by Lemma \ref{condpar}.  Since
$\proj_PR\supseteq\proj_PQ=P$, we have  $\proj_PR=P$ and hence $P'$ is parallel to $P$ by Proposition \ref{Sa46}.

We already have $\proj_{R(P_{i-1},P_i)}P_0=P_{i-1}$ for all $1\leq i\leq k$ by hypothesis. For all $1\leq i\leq
l$, we have $\proj_{R(T_{i-1},T_{i})}P=\proj_{R(T_{i-1},T_{i})}\proj_RP=T_{i-1}$ by Proposition \ref{Sa48}
and because $T_0,T_1,\ldots,T_l$ is a compatible path. This concludes the proof.
\end{Bw}

\begin{Le} \label{le515n}
Let $R$ be a rank $2$ residue and let $P, Q$ be two parallel panels contained in $R$.
Then either $P=Q$ or $R$ is spherical and  $P$ and $Q$ are opposite in $R$.
\end{Le}

\begin{Bw}
Let $J$ be the type of $R$ and let $P$ be an $s$-panel and let $Q$ be a $t$-panel
Then $w:= \delta(P,Q) \in X_s \cap W_J$.
If $l(w) =0$, then $w  = 1_W$ and $P=Q$. Suppose that
$w \neq 1_W$ and let $u \in S$ be such that $l(wu) = l(w)-1$.
As $w \in W_J$ it follows that $u \in J$; moreover $u \neq t$ because
$l(wt) = l(w)+1$ which implies that $J = \{ t,u \}$.
By Lemma \ref{le512n} it follows that $J \subset S^-(wt)$. By Lemma \ref{le53n}
it follows that $J$ is a spherical subset of $S$ and, as $wt \in W_J$,
that $wt = r_J$.
Now $r_J s r_J = r_J^{-1}sr_J = (wt)^{-1}swt = tw^{-1} swt = t^3 = t$
and therefore the panels $P$ and $Q$ are of opposite type with respect to $J$.

Let $c \in P$, put $d := \proj_Qc$ and choose $d' \in Q$ distinct from $d$.
Then $\delta(d,d') = t$ and $\delta(c,d') = \delta(c,d) \delta(d,d') = wt = r_J$.
Hence $c$ and $d'$ are chambers which are opposite inside the residue $R$ and,
as $c \in P$ and $d' \in Q$, it follows that the panels $P$ and $Q$
are opposite inside $R$.
\end{Bw}  

%Indeed, since $P$ and $Q$ are parallel, it follows that $c' := \proj_Pd'$ and $c$ are distinct.
%As $l(c,d') = l(c,d)+1 = l(c',d')+1 = l(c',d)$, it follows that $c'$ and $d$ are both on distinct
%minimal galleries joining $c$ and $d'$. Hence there are two minimal galleries joining $c$ and $d'$
%which shows that $c$ and $d'$ are opposite in $R$.

%Since $c$ was chosen arbitrarily, it follows that for each chamber $c$ in $P$ there is
%a chamber $d'$ in $Q$ which is opposite $c$ inside $R$ and hence $P$ and $Q$
%are opposite panels in $R$.
%\end{Bw}

\begin{Le}\label{comppath}
Two panels are parallel if and only if there exists a compatible path in $\Gamma$ from one to the other.
\end{Le}
\begin{Bw}
The right to left implication will be proved by an induction on the length of the path. If the path has length
one, the result is obvious since opposite panels in a residue are parallel. Assume we have proved the result
for all paths of length strictly less than $k$, and assume $P=P_0,P_1,\ldots P_k=Q$ is a compatible path in
$\Gamma$. By induction $P$ is parallel to $P_{k-1}$. We have $\proj_Q=\proj_Q\proj_{R(P_{k-1},P_k)}$ by Proposition \ref{Sa48}, and so
$\proj_QP=\proj_QP_{k-1}$ which is equal to $Q$ since $P_{k-1}$ and $Q$ are parallel. By Lemma \ref{condpar},
that means $P$ and $Q$ are parallel.

The left to right implication will  be proved by an induction on the numerical distance between the two panels. Let $P,Q$
be two parallel panels. If $l(\delta(P,Q))=0$ then $P=Q$ and the trivial path $P=P_0=Q$ is compatible. Suppose
$l(\delta(P,Q))=l>0$ and the result is proved for all parallel panels at distance strictly less than $l$.
Choose $c\in P$ and let $d=\proj_Qc$. There exists a chamber $e$ adjacent to $d$ such that  $l(c,d)=
l(c,e)+1$. Let $R$ be the unique rank 2 residue containing $Q$ and $e$. By Lemma \ref{projpanels},
$\proj_RP=Q'$ is a panel parallel to $P$ and to $Q$. Since there is a chamber in $R$ closer to $P$ than $d$,
$Q$ cannot be equal to $Q'$ and so they are opposite in $R$ by the previous lemma. Moreover
$l(\delta(P,Q')) <l(\delta(P,Q))$. By induction, there exists a compatible path  $P=P_0,P_1,\ldots P_k=Q'$.
Since $R=R(Q',Q)$, the path $P=P_0,P_1,\ldots P_k,Q$ is compatible.
\end{Bw}

\begin{Le} \label{le517n}
Let $P,Q$ be parallel panels of type $s$ and $t$, respectively, and let $u \in S \setminus \{ t \}$.
Then the following are equivalent:
\begin{itemize}
\item[a)] $u \in S^-(\delta(P,Q))$;
\item[b)] There exists a compatible path $P=P_0,\ldots,P_k=Q$
from $P$ to $Q$ such that $R(P_{k-1},Q)$ has type $\{t,u\}$.
\end{itemize}
Moreover, if this is the case, then $\delta(P,P_{k-1}) = wtr_{\{u,t\}}$
and in particular $l(\delta(P,P_{k-1})) < l(\delta(P,Q))$.
\end{Le}

\begin{Bw}
Let $w = \delta(P,Q)$ and $c \in P$, $J:= \{ u,t \}$ and 
$R$ the $J$-residue containing $Q$. We put $d:= \proj_Qc$, $e:= \proj_Rc$ and $T:=\proj_RP$.
Note that $T$ is parallel to $P$ and to $Q$ by Lemma \ref{projpanels}. As $e \in T$ we have $e = \proj_T e = \proj_T \proj_Rc
=\proj_Tc$ by Proposition \ref{Sa46}. Note that  Proposition \ref{Sa46}  also implies
$d = \proj_Q \proj_R c = \proj_Q e$.

We first show that Assertion a) implies Assertion b). 
Suppose $l(wu)=l(w)-1$, let $U$ be the $u$-panel containing $d$
and $x:= \proj_Uc$. Then $x \neq d$ because $\delta(c,d) = w$ and $l(wu) = l(w)-1$.
It follows that $\proj_Ue = \proj_U\proj_Rc = \proj_Uc = x \neq d$ and in particular
$e \neq d$. As $d = \proj_Qe$, it follows that $\delta(T,Q) \neq 1_W$.
Hence $T$ and $Q$ are opposite in $R$ by Lemma \ref{le515n}.
As $T$ is parallel to $P$, Lemma \ref{comppath} implies that there exists a compatible path $P = P_0,\dots,P_l=T$; setting
$l=k+1$ and
$P_k:= Q$ yields Assertion b).

We now prove that  Assertion b) implies  Assertion a) and the remaining assertions.
Suppose that there exists a compatible path $P=P_0,\ldots,P_k =Q$
such that $R=R(P_{k-1},P_k)$. Then $T=P_{k-1}$.
As $Q$ is opposite $T$ in $R$, we have $\delta(T,Q) = r_Jt$.
Note that $l(r_Jtu) = l(r_Jt)-1$ and that $l(r_Jt) \geq 1$.
We also have that $l(c,d) = l(c,e) + l(e,d)$ and
$w = \delta(P,Q) = \delta(c,d) = \delta(c,e) \delta(e,d) = w_J \delta(T,Q) = w_Jr_Jt$.
We recall that $l(w_Jw') = l(w_J) + l(w')$ for all $w' \in W_J$.
Hence we have
$l(wu) = l(w_Jr_Jtu) = l(w_J) + l(r_Jtu) = l(w_J) + l(r_Jt) - 1 = l(w_Jr_Jt)-1=l(w)-1$ which
yields Assertion a).
We also have $\delta(P,T) = \delta(c,e) = w_J = w_J r_J t t r_J = wtr_J$
and $l(\delta(P,T)) = l(w_J) <l(w_J) + l(r_Jt) =  l(w_Jr_Jt) = l(\delta(P,Q))$.
This finishes the proof.  
\end{Bw}

\begin{Le} \label{le518n}
Let $P$ and $Q$ be parallel panels such that $P \cap Q \neq \emptyset$. Then $P=P_0=Q$ is the only compatible path from
$P$ to $Q$.
\end{Le}
\begin{Bw}
Let $s$ be the type of $P$ and
let $c \in P \cap Q$. Then $\proj_Qc =c$ and therefore $\delta(P,Q) = 1_W$. It follows from
Lemma \ref{parallelpanels} that $Q$ is also of type $s$ and therefore $P=Q$.
Let $P=P_0,\ldots,P_k =Q$ be a compatible path and suppose that $k \geq 1$.
Put $R := R(P_{k-1},P_k)$. By the compatibility of the path we know that
$\proj_RP = P_{k-1}$ and that $P_{k-1}$ is opposite $P_k=Q$ in $R$.
As $P=Q$ is contained in $R$ it follows that $P=\proj_RP$, and so $P=Q=P_{k-1}$. This is a contradiction since a panel cannot be opposite to itself in a rank 2 residue. 
We conclude that $k=0$.
\end{Bw}
 
\begin{Le} \label{le519n}
Let $P$ and $Q$ be parallel panels of type $s$ and $t$ respectively.
Let $c \in P$ and put $d:= \proj_Qc$ and let $E_1(d)$ be the union
of all panels containing $d$. If $l(c,e) \leq l(c,d)$, for all $e \in E_1(d) \setminus Q$,
then ${\cal B}$ is spherical and $P$ and $Q$ are opposite panels.
\end{Le}
\begin{Bw}
Let $w := \delta(P,Q) = \delta(c,d)$. Let $u \in S$ be distinct from $t$
and let $U$ be the $u$-panel containing $d$. We claim that
$x:= \proj_Uc \neq d$. Indeed, if $x=d$ and $y \neq d$ is a chamber
in $U$, then $l(c,y) = l(c,d)+1$ and $y  \in E_1(d) \setminus Q$ which
contradicts our assumption. Hence $x \neq d$ and therefore $u \in S^-(w)$.
It follows now from Lemma \ref{le512n} that $S \subset S^-(wt)$. Hence
$(W,S)$ is a spherical Coxeter system  and $wt$ is the longest element
in $W$, by Lemma \ref{le53n}. Thus  ${\cal B}$ is a  spherical building.

Let $d' \neq d$ be a chamber in $Q$. Then $\delta(c,d') = wt$
and therefore $d$ and $d'$ are opposite chambers in ${\cal B}$.
Moreover, since $wt$ is an involution, we have $(wt)t(wt) = wtttw^{-1} = s$
and therefore $P$ and $Q$ are of opposite types. Hence they are opposite.
\end{Bw}

%\begin{Le}
%Let $R$ be a rank 3 residue and let $P$ and $Q$ be two parallel panels.
%Suppose that  there two compatible paths $P = P_0,\ldots,P_k = Q$
%and $P = Q_0,\ldots, Q_l = Q$ with $k,l \geq 1$ such that $T:= R(P_{k-1},P_k)
%\neq R(Q_{l-1},Q_l) =:T'$. Then $l = k$ and $P_i$ is opposite $Q_{k-i}$
%for $0 \leq i \leq k$.
%\end{Le}
%\begin{Bw}
%Let $c \neq c'$ be chambers in $P$ and put $d:= \proj_Qc, d':= \proj_{Q'}c'$.
%Let $Q'$ be a panel containing $d$ and discinct from $Q$.
%Suppose $Q'$ is contained in $T$. Then $T'

\begin{Le}\label{comprk3}
Let $R$ be a spherical rank 3 residue in ${\cal B}$ and let $P$, $Q$ be two parallel panels in $R$. If there is
more than one compatible path contained in $R$ from $P$ to $Q$, then $P$ and $Q$ are opposite in $R$ and there
are exactly two such paths. Moreover these two paths have the same length.
\end{Le}
\begin{Bw}
Let $P=P_0,P_1,\ldots P_k=Q$ and  $P=P'_0,P'_1,\ldots P'_l=Q$ be two distinct compatible paths in $R$.
We know by  Lemma \ref{le518n}, that $P \neq Q$ and that $k,l \geq 1$.
Let
$Q'=P_{k-i}=P'_{l-i}$ such that $P_{k-j}=P'_{l-j}$ for all $0\leq j\leq i$ and $P_{k-i-1}\neq P'_{l-i-1}$.
Therefore $R(P_{k-i-1},P_{k-i})\neq R(P'_{l-i-1},P'_{l-i})$.

Choose $c\in P$ and let $d=\proj_{Q'}c$. Suppose that $P$ and $Q'$ are not opposite in $R$. It follows from the 
previous lemma  that there
exists a chamber $e$ not in $Q'$ adjacent to $d$ such that $l(c,e)= l(c,d)+1$.
%Let $T \neq Q$ be a panel containing $d$ and suppose that $l(\delta(c,\proj_Tc)) = l(\delta(c,d)-1$  
%Let $d \neq d' \in Q$ and let $x$ be a chamber which is opposite $c$ inside $R$ such that $d'$
%is on a minimal gallery joining $c$ and $x$. If $x = d'$ then it follows that $P$ and $Q$ are opposite
%inside $R$ contradicting our assumption.  
Since there
are only two rank 2 residues in $R$ containing a given panel, the rank 2 residue $R'$ containing $Q'$ and $e$
must be either $R(P_{k-i-1},P_{k-i})$ or $R(P'_{l-i-1},P'_{l-i})$. Without loss of generality we can assume
$R'=R(P_{k-i-1},P_{k-i})$. Then  $\proj_{R'}P=P_{k-i-1}$, which is opposite $Q'$ in $R'$. Let $c'=\proj_{R'}c$. It follows from Part d) of Proposition \ref{sa42} that $l(c',e)= l(c',d)+1$, in contradiction with the fact that $P_{k-i-1}$ is opposite to $Q'$ in $R'$. Therefore $P$ and $Q'$ are opposite. Since $Q'$ cannot be the projection of $P$ on any
rank 2 residue containing it, $Q'$ must be equal to $Q$. Using this and Lemma \ref{comppath} in the building
$R$, we conclude that for two non-opposite parallel panels of $R$, there is exactly one compatible path in $R$
from one to the other.

Let $P$ and $Q$ be opposite in $R$ and let $R'$ be a rank 2 residue in $R$ containing $Q$. Then there is
exactly one compatible path $P=P_0,P_1,\ldots P_k=Q$  such that $R'=R(P_{k-1},P_{k})$. Indeed
$P_{k-1}=\proj_R'P$ is determined and there is only one compatible path between $P$ and $P_{k-1}$ since they
are not opposite. Since there are two rank 2 residues containing $Q$ in $R$, there are exactly two compatible
paths in $R$ from $P$ to $Q$.

A rank 3 spherical residue of a thick building is of type $A_3$, $C_3$, $A_1\oplus A_1\oplus A_1$ or $A_1\oplus I_n$.
Knowing the distance between two opposite panels in $R$ and in all rank 2 residues of $R$, it is easy to
determine the length of compatible paths between opposite panels and see that the two compatible paths have the same length. That length is given in the table below.

\begin{center}
\begin{tabular}{c|ccc|c| p{6cm}}
type $J$&$o(s_1s_2)$&$o(s_2s_3)$&$o(s_1s_3)$&$l(r_J)$&length of compatible paths between opposite panels\\
\hline
$A_3$&$3$&$3$&$2$& $6$&$3$\\
$C_3$&$3$&$4$&$2$& $9$&$4$\\
$A_1\oplus A_1\oplus A_1$&$2$&$2$&$2$&$3$ &$2$\\
$A_1\oplus I_n$&$2$&$n$&$2$&$n+1$ &  $n$ for $s_1$-panels\\
&&&& &  $2$ for $s_2$-panels and $s_3$-panels\\
\hline
\end{tabular}
\end{center}

\end{Bw}

\begin{Le} \label{le521n}
Let $P,Q$ be parallel panels contained in a common residue $R$ and let
$P=P_0,P_1,\ldots,P_k=Q$ be a compatible path. Then $P_i \subset R$
for all $0 \leq i \leq k$.
\end{Le}
\begin{Bw}
We use induction on $k$ and remark that the assertion trivially holds
for $k \leq 1$. Suppose $k \geq 2$ and put $T:=P_{k-1}$ and $R' := R(T,Q)$.
Then $R' \cap R$ is a residue or rank at most 2 containing $Q$.
If it is strictly smaller than 2 we have $Q=R' \cap R$ and $\proj_{R'} R = Q$.
Since $P \subset R$ it follows that $\proj_{R'}P=Q \neq T$, which contradicts
the compatibility of the path. We conclude that $R' \subset R$, which implies
$T \subset R$. Applying induction to the compatible path
$P=P_0,\ldots,P_{k-1} =T$ yields the claim.
\end{Bw}

\begin{Le}\label{compatible}
Let $P$, $Q$ be two parallel panels of ${\cal B}$. Then all compatible paths from $P$ to $Q$ have the same
length.
\end{Le}
\begin{Bw}
Let $P$ be an $s$-panel and let $Q$ be a $t$-panel and let $w:= \delta(P,Q)$.
We will prove the lemma by induction on $l(w)$.

If $l(w)=0$, then $P=Q$ and the trivial path $P=P_0=Q$ is the only compatible path from $P$ to $Q$, by Lemma \ref{le518n}.
Assume $l(w)=L>0$ and we have proved the result for all parallel panels at distance strictly less
than $L$. Take two compatible paths from $P$ to $Q$: $P=P_0,P_1,\ldots, P_k=Q$ and  $P=P'_0,P'_1,\ldots
,P'_l=Q$. If $P_{k-1}=P'_{l-1}=Q'$, then $l(\delta(P,Q'))<L$ and so $k-1=l-1$ and we are done.

Assume now  $P_{k-1}\neq P'_{l-1}$, so that $R_1:= R(P_{k-1},P_{k})\neq R(P'_{l-1},P'_{l})=:R_2$, and let $R$ be the rank
3 residue containing these two rank 2 residues. As $Q \subset R_i$ the residue $R_i$ has type $\{t,u_i\}$ for $i=1,2$;
thus we have
$u_1 \neq u_2 \in S$ and $\{t,u_1,u_2\}$ is the type of $R$. 
By Lemma \ref{le517n} it follows that $\{ u_1,u_2 \} \subset S^-(w)$ and therefore,
by Lemma \ref{le512n}, the residue $R$ is spherical.

Let $Q'$ be the projection of $P$ on $R$. By Lemma
\ref{projpanels}, $Q'$ is parallel to $P$ and $Q$. 
As there exists a compatible path from $P$ to $P_{k-1}$, the panel $P_{k-1}$ is parallel to $P$ by Lemma \ref{comppath}
and as $P_{k-1} \subset R$, it is also parallel to $Q'$ by Lemma \ref{projpanels}. 
Hence there exists a compatible path  $Q'=T_0,T_1,\ldots,T_m=P_{k-1}$
from $Q'$ to $P_{k-1}$ and all the $T_i$ are contained in $R$ by the previous lemma.
Similarly we have a compatible path $Q'=T'_0,T'_1,\ldots,T'_n=P'_{l-1}$ in $R$.
We have $P_{k-1}=\proj_{R(P_{k-1},P_{k})}P=\proj_{R(P_{k-1},P_{k})}\proj_RP=\proj_{R(P_{k-1},P_{k})}Q'$, and so
$Q'=T_0,T_1,\ldots,T_m,P_k=Q$ is a compatible path. By similar arguments,  $Q'=T'_0,T'_1,\ldots,T'_n,P'_l=Q$ is
also a compatible path. As $P_{k-1} \neq P'_{l-1}$, these paths are distinct. 
It follows from Lemma \ref{comprk3} that $Q'$ and $Q$ are opposite panels in $R$
and that these two paths in $R$  have the same length, hence $m=n$.

By Lemma \ref{comppath}, there is a compatible path from $P$ to $Q'$, denoted by $P=S_0,S_1,\ldots,S_j=Q'$. By
Lemma \ref{projpanels}, we have, on the one hand, that the paths $P=S_0,S_1,\ldots,S_j=Q'=T_0,T_1,\ldots,T_m=P_{k-1}$ and
$P=S_0,S_1,\ldots,S_j=Q'=T'_0,T'_1,\ldots,T'_m=P'_{l-1}$ are both compatible of length $j+m$. On the other hand
$P=P_0,P_1,\ldots, P_{k-1}$ and $P=P'_0,P'_1,\ldots ,P'_{l-1}$ are also compatible paths. Since
$l(\delta(P,P_{k-1}))<L$ and $l(\delta(P,P'_{l-1}))<L$, we can use the hypothesis of induction, and so
$k-1=j+m$ and $l-1=j+m$. We conclude that $k=l$.
\end{Bw}

\begin{De}
By Lemma \ref{compatible}, we can define the \emph{compatible distance between
two parallel panels $P$ and $Q$} as the length of a compatible path joining them. It will
be denoted by $L(P,Q)$. %Note, that the previous lemma shows that the compatible distance between two parallel panels does not depend on the choice of a compatible path.
\end{De}

\begin{Le}\label{compatiblew} Let $w \in X_s$ and let $P,P'$ be $s$-panels and
$Q,Q'$ be $w^{-1} s w$-panels such that
$\delta(P,Q) = w = \delta(P',Q')$. Then $L(P,Q) = L(P',Q')$.
Moreover, if $(J_1,\ldots,J_k)$ is  the type of a compatible path
from $P$ to $Q$, then there exists a compatible path from $P'$ to $Q'$
of the same type.
\end{Le}
\begin{Bw} This follows by induction on $l(w)$ using Lemma \ref{le517n}.
\end{Bw}

\begin{De}
Let $w \in X_s$. Then we define its \emph{$s$-compatible
length}, denoted by $L_s(w)$, as the compatible distance
between an $s$-panel $P$ and an $w^{-1} s w$-panel $Q$
such that $\delta(P,Q)=w$.
\end{De}

\begin{Sa} \label{Sa526n}
Let $s \in S$ and  $w_1,w \in X_s$ such that  $w_1\prec w$. Put $w_2 := w_1^{-1}w, u:= w_1^{-1}sw_1, t:=w^{-1}sw$. 
Let $c,e,d$ be chambers such that $\delta(c,e) = w_1$ and $\delta(e,d) = w_2$
and let $P$ be the $s$-panel containing $c$, $U$ the $u$-panel containing $e$ and $Q$ be the $t$-panel
containing $d$. Then the following holds.
\begin{itemize}
\item[a)] $w_2 \in X_u$;
\item[b)] the three panels are pairwise parallel and we have $\delta(P,U)= w_1,\delta(U,Q)=w_2$
and $\delta(P,Q) = w$;
\item[c)] $\proj_Q^P = \proj_Q^U \circ \proj_U^P$.
\end{itemize}
\end{Sa}
\begin{Bw} 
First note that $u$ and $t$ are in $S$, since $w_1,w \in X_s$. Since $w_1\prec w$, we also have that $l(w)=l(w_1)+l(w_2)$.
We easily see that $uw_2=w_2t$, so $ w_2^{-1}uw_2\in S$. Since $w \in X_s$ we also have $l(wt)=l(sw)=l(w)+1$.  
Therefore $ l(w_1)+l(w_2)+1=l(w_1w_2t)\leq l(w_1)+l(w_2t)\leq l(w_1)+l(w_2)+1 $, and so $l(w_2t)=l(w_2)+1$. Hence Part a) holds.
 Part b) follows from Part a) and the second assertion of Lemma \ref{parallelpanels}.
Let $x \in P$. Then $\delta(x,\proj_Ux) = w_1$ and $\delta(\proj_Ux,\proj_Q \proj_Ux) = w_2$,
by Part b). As $l(w_1w_2) = l(w_1)+l(w_2)$ it follows that
$\delta(x,\proj_Q \proj_Ux) = w_1w_2 = w$. Now $\proj_Qx$ is the unique chamber $y$ in $Q$
such that $\delta(x,y) = \delta(P,Q) = w$ and therefore $\proj_Qx = \proj_Q \proj_Ux$. Hence Part c) holds.
\end{Bw}

%%%%%%%%%%%�������������������

\section{Projectivities between panels}

\label{bijections}

Throughout this section, ${\cal B} =({\cal C},\delta)$ is  a building of type $(W,S)$ and $f$ is a codistance
on ${\cal B}$. Moreover, it is always assumed that ${\cal B}$ satisfies the Conditions~(3-sph), 
(lco) and (lsco) of Theorem \ref{th52}. By the latter result we have in particular that $f^{\op}$ is 
simply connected.

We first recall some facts about codistances and fix further notation.

Let $c \in {\cal C}$ and $w \in W$ be such that $l(f(c)w) = l(f(c)) + l(w)$. By Lemma \ref{uniquec}
there is a unique chamber $d \in {\cal C}$ such that $\delta(c,d) = w$ and $f(d) = f(c)w$.
We denote this chamber by $\pi(c,w)$. Note that $\pi(c,w)$ is  defined for all $w$
if $c \in f^{\op}$.

The following observation is immediate.

\begin{Le} \label{le61n}
Let $c \in f^{\op}$,  $w_1 \prec w \in W$ and put $w_2 := w_1^{-1}w$.
Then $\pi(\pi(c,w_1),w_2) = \pi(c,w)$.
\end{Le}

\begin{De}
We will say that a residue $R$ is \textit{in} $f^{\op}$, or in $f^{\op}_c$, respectively, if it contains a chamber in
$f^{\op}$, or in $f^{\op}_c$, respectively. For $s\in S$, let ${\cal P}^{\op}_s(f)$ and ${\cal P}^{\op}_{s,c}(f)$, respectively, be
the set of all $s$-panels in $f^{\op}$ and $f^{\op}_c$, respectively.
\end{De}

Notice that all chambers of a panel $P$  in $f^{\op}$ are in $f^{\op}$ except for one, namely $\proj_Pf$.

\begin{Sa}\label{pi}
Let $P\in {\cal P}^{\op}_s(f)$, $w\in X_s$ and $t=w^{-1}sw$. Let $P'$ be a $t$-panel  with  $\delta(P,P')=w$.
Then the following conditions are equivalent:
\begin{itemize}
\item[a)] $P'$ contains a chamber with $f$-value $w$;
\item[b)] $f(x)\in\{w,wt\}$ for $x\in P'$ and exactly one chamber of $P'$ has $f$-value $wt$;
\item[c)] $P\in  {\cal P}^{\op}_{s,x}(f)$ for all chambers $x$ of $P'$;
\item[d)] $P\in  {\cal P}^{\op}_{s,x}(f)$ for some chamber $x$ of $P'$;
\end{itemize}
There exists exactly one panel $P'$ satisfying these conditions.
\end{Sa}
\begin{Bw}
Conditions a) and b) are equivalent by the definition of a codistance. Assume $P'$ satisfies b). Let $x$ be a
chamber with $f$-value $w$ in $P'$. Then $\delta(\proj_Px,x)=w=f(x)$. Since $\proj_Px$ cannot be equal to
$\proj_Pf$ (otherwise the chamber in $P'$ with $f$-value $wt$ would be at distance $l(w)$ from a chamber in $f^{\op}$
yielding a contradiction in view of Lemma \ref{galltofop}),
$\proj_Px\in f^{\op}_x$ and $P\in  {\cal P}^{\op}_{s,x}(f)$. Now let $z=\proj_{P'}f$ be the unique chamber in $P'$ with
$f$-value $wt$. If $y$ is any chamber of $P$ in  $f_{\op}$, $\delta(y,z)=wt$, so $y\in f^{\op}_z$ and $P\in
{\cal P}^{\op}_{s,z}(f)$. Obviously c) implies d). Now assume $P'$ satisfies d). Then $P$ contains $y\in
f^{\op}$ and $\delta(y,x)=f(x)$. If $y=\proj_Px$, then $\delta(y,x)=w$; if $y \neq \proj_Px$, then $\delta(y,x)=sw=wt$. In both
cases, $P'$ contains a chamber with $f$-value $w$.

We prove the existence of such a panel. Let $p$ be the unique chamber of $P$ not in $f^{\op}$ (so
$f(p)=s$). As $l(f(p)w)=l(sw)=l(w)+1$, Lemma~\ref{uniquec} implies the existence of a unique chamber $c$ with
$sf(c)=w=\delta(p,c)$. Let $P'$ be the $t$-panel on $c$. By Lemma \ref{parallelpanels}, $P'$ is parallel to $P$
and $\delta(P,P')=w$. It obviously satisfies a).

Now we want to show that $P'$ is unique. Let $Q$ be a $t$-panel with  $\delta(P,Q)=w$ satisfying b). Let $x$ be
the chamber of $Q$ with $f$-value $wt=sw=f(p)w$. We have $\proj_Px=p$, so $\delta(p,x)=w$. By Lemma
\ref{uniquec}, a chamber with that property is unique.
 Therefore $x=c$ and $Q=P'$.
\end{Bw}

\begin{De} For  $P\in {\cal P}^{\op}_s(f)$, $w\in X_s$ and $t=w^{-1}sw$, we denote the unique $t$-panel $P'$ with $\delta(P,P')=w$ satisfying the equivalent conditions a) up to d) of Proposition~\ref{pi} by $\pi(P,w)$.
\end{De}

\begin{Le} \label{projlemma}
Let $P\in {\cal P}^{\op}_s(f)$, $w\in X_s$ and $t=w^{-1}sw$ and $P' := \pi(P,w)$.
Then $\proj_P \proj_{P'} f = \proj_P f$.
\end{Le}
\begin{Bw}
As $l(\proj_{P'} f) = l(w)+1$ and $l(\proj_P \proj_{P'} f,  \proj_{P'} f) = l(w)$, it follows that
$l(f(\proj_P \proj_{P'} f)) \geq 1$ and hence $f(\proj_P \proj_{P'} f) =s$. The claim follows
because $\proj_Pf$ is the unique chamber in $P$ having $f$-value $s$.
\end{Bw}

\begin{Le}\label{revpi}
Let $Q$ be a $t$-panel of ${\cal B}$ and let $w$ be the shortest word of $\{f(x)|x\in Q\}$. Suppose
$wtw^{-1}:=s\in S$. Then there exists an $s$-panel $P\in {\cal P}^{\op}_s(f)$ such that $Q=\pi(P,w)$.
\end{Le}
\begin{Bw}
Since $w^{-1}sw=t$ and $l(sw)=l(wt)=l(w)+1$, we have $w\in X_s$. Let $x$ be a chamber of $Q$ with $f(x)=w$. Let
$y\in f^{\op}_x$ so that $\delta(y,x)=w=f(x)$. Let $P$ be the $s$-panel on $y$. By construction $P$ is a panel
in ${\cal P}^{\op}_{s,x}(f)$ which is parallel to $Q$ by Lemma \ref{parallelpanels}. Moreover $\delta(P,Q)=w$,
hence by Proposition \ref{pi}, $Q=\pi(P,w)$.
\end{Bw}

\begin{Le} \label{le65n}
Let $c \in f^{\op}$, $s \in S$, $w \in X_s$ and put $t := w^{-1}sw$.
Let $P$ be the $s$-panel containing $c$ and let $Q$ be the $t$-panel containing
$d:=\pi(c,w)$. Then $Q = \pi(P,w)$ and $d=\proj_Qc$. Moreover,
$\proj_P \proj_Qf= \proj_Pf$ and $\proj_Q \proj_P f = \proj_Q f$.
\end{Le}
\begin{Bw}
By its description, $Q$ is a $t$-panel containing a chamber $d$ such that $\delta(c,d) = w = f(d)$.
Since $w \in X_s$, it follows that $Q$ is parallel to $P$ and $\delta(P,Q) = w$; hence $Q = \pi(P,w)$
by Condition a) of Proposition \ref{pi}. As $\delta(c,d) = w = \delta(P,Q)$ we have
$d = \proj_Qc$.

Let $x \in Q$. Then, by Condition b) of Proposition \ref{pi}, we have $f(x) \in \{ w,wt \}$.
As $l(wt) = l(w)+1$, we have $x = \proj_Q f$ if and only if $f(x) = wt$.
Let $q:=\proj_Qf$ and $p := \proj_Pq$. Then $\delta(p,q) = w$ and therefore
$l(\delta(p,q)) < l(f(q))$. It follows that $p$ is not in $f^{\op}$ by Lemma \ref{galltofop}.
As $f(y) \in \{ 1_W,s \}$ for all $y \in P$, we conclude that $f(p) =s$ and therefore
$p = \proj_Pf$. The second equality follows from the fact that 
$\proj_Q^P$ and $\proj_P^Q$ are inverse to each other.
\end{Bw}

%\begin{Sa}
%Let $s \in S$,  $w_1 \prec w \in X_s$ and $P,P' \in {\cal P}_s^{\op}(f)$ such that $P 

\begin{De}
For $P,Q\in {\cal P}^{\op}_s(f)$ and $w\in X_s$, we write $P\equiv_w Q$ if $\pi(P,w)=\pi(Q,w)$. This
is an equivalence relation on ${\cal P}^{\op}_s(f)$. For $P\equiv_w Q$, we put $\beta(P,Q,w)$ the bijection
from $P$ to $Q$ defined by $\proj_Q\proj_{\pi(P,w)}$.
\end{De}

Notice that $\beta(Q,P,w)\beta(P,Q,w)=1_P$ and that, by Lemma \ref{projlemma}, $\beta(P,Q,w)$ maps $\proj_Pf$ onto
$\proj_Qf$ via  $\proj_{\pi(P,w)}f$.

\begin{Sa} \label{extension}
Let $s \in S$, $w_w,w\in X_s$, and suppose $w_1 \prec w$. Let $P,P' \in {\cal P}_s^{\op}(f)$ such that $P \equiv_{w_1} P'$.
Then $P \equiv_w P'$ and $\beta(P,P',w) = \beta(P,P',w_1)$.
\end{Sa}
\begin{Bw}
Let $U = \pi(P,w_1) = \pi(P',w_1)$, let $Q = \pi(P,w)$ and let $c \in P \cap f^{\op}$.
Put $t := w^{-1}sw, u:= w_1^{-1}sw_1, w_2:= w_1^{-1}, e:=\pi(c,w_1)$ and $d:=\pi(c,w)$.
By Lemma \ref{le65n} $U$ is the $u$-panel containing $e$ and $Q$ is the $t$-panel containing $d$. By Lemma \ref{le61n}, we have $\pi(e,w_2) = d$, and in particular $\delta(e,d)= w_2$.
It now follows from Proposition \ref{Sa526n} that $U$ is parallel to $Q$ and that
$\proj_Q^P = \proj_Q^U \circ \proj_U^P$.

Put $c':=\proj_{P'} e$. As $\delta(c',e) = \delta(P',U) = w_1 = f(e)$,
we have $c' \in f^{\op}$ and $e = \pi(c',w_1)$. Now $\pi(c',w) = \pi(\pi(c',w_1),w_2) = \pi(e,w_2) = d$ 
which implies that $Q = \pi(P',w)$ and shows that $P \equiv_w P'$.  We now apply 
Proposition~\ref{Sa526n} again to see that
$\proj_Q^{P'} = \proj_Q^U \circ \proj_U^{P'}$.
As $\proj_Q^U$ and $\proj_U^Q$ are mutually inverse bijections, it follows that
$\beta(P,P',w) = \beta(P,P',w_1)$.
\end{Bw}

\begin{Le} \label{le69nn}
Let $c \in f^{\op}$, $w \in W$ and let $R$ be a spherical $J$-residue containing $\pi(c,w)$.
Then $\proj_Rf = \pi(c,w^J)$ and $\proj_Rc = \pi(c,w_J)$.
\end{Le}
\begin{Bw}
Let $d:=\pi(c,w)$. As $d \in R$ we have $f(R)=f(d)W_J=wW_J$ by Lemma~\ref{lem:imfR}. It follows
that $\delta(c,\proj_Rc) = w_J$ and $\delta(\proj_Rc,d)=w_J^{-1}w$ and that $\proj_Rc$ is on a minimal
gallery joining $c$ and $d$. As $d = \pi(c,w)$, it follows that $\proj_Rc=\pi(c,w_J)$.
Now $\pi(c,w^J) = \pi(c,w_Jr_J)$ and, as $w_J \prec w^J$, we have
$ \delta(\pi(c,w_J),\pi(c,w^J)) = r_J$. Consequently, $\pi(c,w^J) \in R$ and $\pi(c,w^J)$ is the 
unique element in $R$ having $f$-value $w^J$, which yields $\pi(c,w^J) = \proj_Rf$.
\end{Bw}

\begin{Le} \label{le610nn}
Let $P \in {\cal P}^{\op}_s(f)$, $w \in X_s$, and let $R$ be a spherical $J$-residue containing $\pi(P,w)$.
Then $w_J, sw^J \in X_s$, $\pi(P,w_J) = \proj_R P$,  and $w_J \prec w \prec sw^J$.
Also, $\pi(P,sw^J)$ is the $t$-panel containing $\proj_Rf$ where $t = (w^J)^{-1} s w^J$.
\end{Le}

\begin{Bw}
We put $Q:= \pi(P,w)$. Let $c \in P \cap f^{\op}$ and put $T_1 := \proj_RP$ and $c_1 := \proj_Rc \in T_1$.
First note  that $T_1$ is parallel to both $P$ and $Q$ by Lemma \ref{projpanels} and that
$w_J = \delta(c,c_1)$ by Part d) of Proposition~\ref{sa42}. As $c_1 \in T_1 \subset R$, we have
$c_1 = \proj_R c = \proj_{T_1} c$ and therefore $w_J = \delta(c,c_1) = \delta(P,T_1) \in X_s$.
By the previous lemma we know that $c_1 = \pi(c,w_J)$ and as $c_1 \in T_1$ we obtain
$T_1 = \pi(P,w_J)$.

Let $u:= w_J^{-1} s w_J$. As $w_J \in X_s$, it follows that $u \in S$ and hence that $T_1$ is a $u$-panel.
As $T_1 \subset R$, we obtain $u \in J$. We put $t := r_J u r_J\in J$ and recall that $w^J = w_Jr_J$.
This yields $t = (w^J)^{-1} s w^J$ and in particular $sw^J = w^Jt$.
As $t \in J$, we have $l(r_Jt) = l(r_J)-1$ and therefore
$l(sw^J) = l(w^Jt) = l(w_J) + l(w_J^{-1} w^Jt) = l(w_J) + l(r_Jt) = l(w_J) + l(r_J)-1 = l(w_Jr_J) -1 = l(w^J)-1$,
hence $sw^J \in X_s$. 
Let $T_2 := \pi(P,sw^J)$ and put $c_2 := \proj_{T_2}c$. Then $T_2$ is a $t$-panel and $c_2 = \pi(c,sw^J)$
by Lemma \ref{le65n}. As $\delta(c_1,c_2) = \delta(\pi(c,w_J),\pi(c, w^Jt) = r_Jt \in W_J$,
it follows that $c_2 \in R$ and therefore $T_2$ is contained in $R$. 
Now $T_2$ is a $t$-panel contained in $R$ which contains a chamber having $f$-value $w^Jt$, 
hence
it contains also the unique chamber in $R$ having $f$-value $w^J$ which is in fact the projection 
of $f$ onto $R$.

It remains to show that $w \prec sw^J$. As $v:= w_J^{-1}w \in W_J$, we have 
$t' := w^{-1}sw = v^{-1}uv \in S \cap W_J = J$ and therefore $sw = wt' \prec w^J$.
As $ l(s(sw)) = l(sw)-1$ and $l(sw^J) = l(w^J)-1$ the assertion follows.
\end{Bw}

\begin{Ko} \label{ko611}
Let $P,Q\in {\cal P}^{\op}_s(f)$, $w \in X_s$ and let $J \subset S$ be spherical.
If $\pi(P,w)$ and $\pi(Q,w)$ are contained in the same $J$-residue $R$,
then $P \equiv_{sw^J} Q$.
\end{Ko}
\begin{Bw}
Put $t := (w^J)^{-1} sw^J$ and let $T$ be the $t$-panel containing $\proj_Rf$.
Then we have $\pi(P,sw^J) = T = \pi(Q,sw^J)$ by the previous lemma.
\end{Bw}

\begin{De} For $P,Q\in {\cal P}^{\op}_s(f)$, we say that $P$ and $Q$ are \emph{$t$-adjacent}, denoted by $P\sim_t Q$, if there exist
$p\in P\cap f^{\op}$ and $q\in Q\cap f^{\op}$ with $p\sim_t q$. Let $P\sim_t Q$, both in ${\cal P}^{\op}_s(f)$.
%, we will define a bijection $\beta(P,Q)$ from $P$ to $Q$. If $t=s$, then $P=Q$ and we take $\beta(f,P,Q)=1_P$.

Let $J=\{s,t\}$ and let $R$ be the $J$-residue containing $P$ and $Q$. As ${\cal B}$ is assumed to be 3-spherical,
the residue  $R$ is spherical and we put $x_J:=sr_J\in X_s$. 
%Notice that $\pi(P,x_J)$ is a panel opposite to $P$ in $R$ which contains a
%chamber with $f$-value $r_J$. By Proposition \ref{projf}, $R$ contains only one such chamber, namely
%$\proj_Rf$. Hence  $\pi(P,x_J)$ is the unique panel of type $x_J^{-1}sx_J=r_J^{-1}sr_J$ containing $\proj_Rf$.
%Similarly  $\pi(Q,x_J)$ is the same panel. 
By the previous corollary we have $P\equiv_{x_J} Q$ and  we put
$\alpha(P,Q):=\beta(P,Q,x_J)$. If $s=t$, $x_J=1_W$, $P=Q=\pi(P,x_J)$ and $\alpha(P,Q)=\id_P$.
\end{De}

Notice that if  $P,Q\in {\cal P}^{\op}_s(f)$ are $t$-adjacent, then $\alpha(P,Q) (\proj_Pf) = \proj_Qf$ in view of Lemma~\ref{projlemma}.

\begin{Le} \label{le612nn}
Let $R$ be a spherical $J$-residue in $f^{\op}$ and let $c,d \in R \cap f^{\op}$.
Let $s \in J$ and let $P$ and $Q$ be the $s$-panels containing $c$ and $d$, respectively.
Let $c=c_0,\ldots,c_m=d$ be a gallery in $R \cap f^{\op}$ and for each
$0 \leq i \leq m$ let $P_i$ be the $s$-panel containing $c_i$.
Then $\alpha(P_{m-1},P_m)\alpha(P_{m-2},P_{m-1})\ldots\alpha(P_1,P_2)\alpha(P_0,P_1) = \beta(P,Q,sr_J)$.
\end{Le}

\begin{Bw}
We put $x_J := sr_J$ and observe that  $T:=\pi(P,x_J) = \pi(P_1,x_J) =\ldots = \pi(P_{m-1},x_J) = \pi(Q,x_J)$ by Corollary \ref{ko611}.

Let $1 \leq i \leq m$. If $c_{i-1}$ is $s$-adjacent to $c_i$, we put $J_i := \{ s \}$; if they are not $s$-adjacent,
then they are $t_i$-adjacent for a unique $t_i \in J$, $t_i\neq s$, and we put $J_i := \{ s,t_i \}$ in this case.
Furthermore, we put $x_i := sr_{J_i} \in X_s$ and observe that $x_i \prec x_J$ for all $1 \leq i \leq m$ by Lemma~\ref{le610nn}. 

We can now apply Proposition \ref{extension} to see that
$\alpha(P_{i-1},P_i) = \beta(P_{i-1},P_i,x_i) = \beta(P_{i-1},P_i,x_J) = \proj_{P_{i-1}}^T \proj_T^{P_i}$,
for all $1 \leq i \leq m$. As $\proj_{P_i}^T$ and $\proj_T^{P_i}$
are inverse bijections for $1 \leq i \leq m-1$, we obtain
 $\alpha(P_{m-1},P_m)\alpha(P_{m-2},P_{m-1})\ldots\alpha(P_1,P_2)\alpha(P_0,P_1) = 
\proj_Q^T \proj_T^P = \beta(P,Q,sr_J)$.
\end{Bw}

\begin{Th}\label{beta}
There exists a unique system of bijections $\beta(P,Q): P \rightarrow Q$ where  $P,Q$ in ${\cal P}^{\op}_s(f)$,
such that the following conditions 
are satisfied for all  $P,Q,R \in {\cal P}^{\op}_s(f)$ :
%For any two $s$-panels $P,Q$ in ${\cal P}^{\op}_s(f)$, we can define a bijection $\beta(P,Q)$ from $P$ to $Q$
%in such a way that the following hold for all $P,Q,R \in {\cal P}^{\op}_s(f)$ :
\begin{itemize}
\item[a)] $\beta(P,P)=1_P$;
\item[b)] $\beta(Q,P)\beta(P,Q)=1_P$;
\item[c)] $\beta(Q,R)\beta(P,Q)=\beta(P,R)$;
\item[d)] $\beta(P,Q)(\proj_Pf)=\proj_Qf$;
\item[e)] if $P$ and $Q$ are $t$-adjacent for some $t \in S$, then $\beta(P,Q) = \alpha(P,Q)$.
\end{itemize}
%We can define an equivalence relation on the chambers in the $s$-panels of ${\cal P}^{\op}_s(f)$ such that each equivalence class meets each $s$-panel of ${\cal P}^{\op}_s(f)$ in exactly one chamber and the set $\{\proj_Pf|P\in{\cal P}^{\op}_s(f)\}$ forms an equivalence class.
\end{Th}
\begin{Bw}
%In order to do this, we will define bijections between panels of  ${\cal P}^{\op}_s(f)$.

%For $P,Q\in {\cal P}^{\op}_s(f)$, we say that they are $t$-adjacent, denoted by $P\sim_t Q$, if there exist
%$p\in P\cap f^{\op}$ and $q\in Q\cap f^{\op}$ with $p\sim_t q$. Let $P\sim_t Q$, both in ${\cal P}^{\op}_s(f)$.
%, we will define a bijection $\beta(P,Q)$ from $P$ to $Q$. If $t=s$, then $P=Q$ and we take $\beta(f,P,Q)=1_P$.
%Let $J=\{s,t\}$ and $R$ be the $J$-residue containing $P$ and $Q$. Since $R$ is spherical by the hypothesis on
%${\cal B}$, $x_J:=sr_J\in X_s$. Notice that $\pi(P,x_J)$ is a panel opposite to $P$ in $R$ which contains a
%chamber with $f$-value $r_J$. By Proposition \ref{projf}, $R$ contains only one such chamber, namely
%$\proj_Rf$. Hence  $\pi(P,x_J)$ is the unique panel of type $x_J^{-1}sx_J=r_J^{-1}sr_J$ containing $\proj_Rf$.
%Similarly  $\pi(Q,x_J)$ is the same panel. Therefore $P\equiv_{x_J} Q$ and  we put
%$\beta(P,Q):=\beta(P,Q,x_J)$. If $s=t$, $x_J=1_W$, $P=Q=\pi(P,x_J)$ and $\beta(P,Q)=1_P$.

%As noticed above, $\beta(P,Q)$ maps $\proj_Pf$ onto $\proj_Qf$.

%We have seen that for, any $x\in R$, $\delta(\proj_Rf,x)=f(\proj_Rf)^{-1}f(x)$ which has length $l(f(\proj_Rf))-l(f(x))$, therefore the distance from $\proj_Rf$ to chambers in $P$ is minimal when $f(x)$ is the longest, that is for $\proj_Pf$. Hence $\proj_Rf$ is projected on $P$ onto $\proj_Pf$ and on $Q$ onto $\proj_Qf$. We conclude that $\beta(P,Q)$ maps $\proj_Pf$ onto $\proj_Qf$.

Let $P,Q\in {\cal P}^{\op}_s(f)$ and choose $p\in P\cap f^{\op}$ and $q\in Q\cap f^{\op}$. By Theorem~\ref{th52}, $f^{\op}$ is connected, and so there exists a gallery $\gamma$ from $p$ to $q$ contained in
$f^{\op}$. Set $\gamma=(x_0=p,x_1,x_2,\ldots,x_n=q)$ and let $X_i$ be the $s$-panel containing $x_i$. By definition
these panels are in ${\cal P}^{\op}_s(f)$ and $X_i\sim_{t_i}X_{i+1}$ for some $t_i\in S$. We define
$\beta(\gamma,P,Q):=\alpha(X_{n-1},Q)\ldots\alpha(X_1,X_2)\alpha(P,X_1)$. By the above comment $\beta(\gamma,P,Q)$
maps $\proj_Pf$ onto $\proj_Qf$. Note also that, if a system of bijections  satisfying the conditions of the theorem exists, then $\beta(P,Q)$ has to coincide
with $\beta(\gamma,P,Q)$ in view of Condition e). This yields already the uniqueness and it remains to show
that the bijection $\beta(P,Q)$ defined by $\beta(\gamma,P,Q)$ does not depend on the choice of $\gamma$.

We will now show that if $\gamma_1$ and $\gamma_2$ are two galleries in $f^{\op}$ from a chamber of $P$ to a
chamber of $Q$, then $\beta(\gamma_1,P,Q)=\beta(\gamma_2,P,Q)$. As  $\beta(P,P)=\id_P$, we can assume that $\gamma_1$
and $\gamma_2$ start and finish with the same chamber. Hence this is equivalent to showing that for a closed
gallery $\gamma$ in $f^{\op}$,  $\beta(\gamma,P,P)=\id_P$.

We recall that the assumptions on the buildings considered in this section allow us to apply
Theorem \ref{th52}; hence $f^{\op}$ is simply 2-connected. Therefore there exists a finite sequence of elementary
homotopies from the closed gallery $\gamma$ to a trivial gallery based in $p\in P$ such that all intermediate
galleries are contained in $f^{\op}$. Since two galleries differing by an elementary homotopy are equal except
in a rank 2 residue, it is enough to show that $\beta(\gamma,P,P)=\id_P$ for a closed galley $\gamma$ in a rank 2 residue in $f^{\op}$.

Let $\gamma=(x_0,x_1,x_2,\ldots,x_n=x_0)$ be a closed gallery in $f^{\op}$ contained in a rank 2 residue $R$ of
type $\{t,u\}$ (where $t$ or $u$ could be equal to $s$). Let the $X_i$'s be defined as above and put $J=\{s,t,u\}$.
In view of Lemma \ref{le612nn} we now have
$\beta(\gamma,X_0,X_0)=\alpha(X_{n-1},X_0)\ldots\alpha(X_1,X_2)\alpha(X_0,X_1)=
\beta(X_0,X_0,sr_J) = \id_{X_0}$.

% Consider two
%consecutive chambers $x_i$ and $x_{i+1}$. They are $v$-adjacent, where $v\in\{t,u\}$. Let $J:=\{s,v\}$ and
%$K=\{s,t,u\}$, which are spherical by the hypothesis on ${\cal B}$. Let $R$ be the $K$-residue containing
%$x_0$. We have  $X_i\sim_{v}X_{i+1}$ and, as seen above $X_i\equiv_{x_J}X_{i+1} $ and
%$\alpha(X_i,X_{i+1})=\beta(X_i,X_{i+1},x_J)$. By Lemma \ref{Xs}, $x_J\prec x_K$, and so, by Proposition
%\ref{extension}, $X_i\equiv_{x_K}X_{i+1}$ and $\alpha(X_i,X_{i+1})=\beta(X_i,X_{i+1},x_K)$. For all
%$j=0,1,\ldots ,n-1$, $\pi(X_j,x_K)$ is the panel of type $x_K^{-1}sx_K$ containing $\proj_Rf$, which is
%opposite to $X_j$ in $R$. This panel wil bel denoted by $\pi(R)$. We have
%$$\begin{array}{ll}
%\beta(\gamma,X_0,X_0)&=\alpha(X_{n-1},X_0)\ldots\alpha(X_1,X_2)\alpha(X_0,X_1)\\
%&=\beta(X_{n-1},X_0,x_K)\ldots\beta(X_1,X_2,x_K)\beta(X_0,X_1,x_K)\\
%&=\proj_{X_0}\proj_{\pi(R)}\ldots\proj_{X_2}\proj_{\pi(R)}\proj_{X_1}\proj_{\pi(R)}\\
%&=\proj_{X_0}\proj_{\pi(R)}=1_{X_0}
%\end{array}$$
%because $\proj_{\pi(R)}\proj_{X_i}=\id_{\pi(R)}$. This completes the proof that
%$\beta(\gamma_1,P,Q)=\beta(\gamma_2,P,Q)$ for $\gamma_1$ and $\gamma_2$  two galleries in $f^{\op}$ from a
%chamber of $P$ to a chamber of $Q$. We put $\beta(P,Q):=\beta(\gamma,P,Q)$ for $\gamma$ any gallery in
%$f^{\op}$ from a chamber of $P$ to a chamber of $Q$.

It is obvious that a), b) and c) are satisfied. Since d) is satisfied for adjacent panels, it will be
satisfied, by induction, for any two panels. Finally, Condition e) is satisfied by the construction of
the system of bijections $\beta(P,Q)$.\end{Bw}

%Let $P\in {\cal P}^{\op}_s(f)$ and $x$ a chamber of $P$.
%We define the equivalence class of $x$ as $\{\beta(P,Q)(x)|Q\in {\cal P}^{\op}_s(f)\}$. It is obvious that this gives an equivalence relation.
%By construction, each class contains one chmaber in each panel of ${\cal P}^{\op}_s(f)$. Moreover, the equivalence class of $\proj_Pf$ is $\{\beta(P,Q)(\proj_Pf|%Q\in {\cal P}^{\op}_s(f)\}=\{\proj_Qf)|Q\in {\cal P}^{\op}_s(f)\}$ by an earlier remark.

\begin{Th}\label{betaw}
Let $(\beta(P,Q))_{P,Q \in {\cal P}^{\op}_s(f)}$ be the unique system of bijections satisfying
the conditions of the previous theorem.
Let $P,P'\in {\cal P}^{\op}_s(f)$  with $P\equiv_{w} P'$ for $w\in X_s$. Then $\beta(P,P')=\beta(P,P',w)$.
\end{Th}
\begin{Bw}
We first consider the special case where $w \in W_J$ for some spherical subset $J$ of $S$ containing $s$.
Let $R$ be the spherical $J$-residue containing $P$ and $P'$. Let $p \in P \cap f^{\op}$ and $p' \in P' \cap f^{\op}$.
As $R \cap f^{\op}$  is connected, there is a gallery $p=x_0,\ldots,x_n=p'$ in 
$R \cap f^{\op}$ and we let $X_i$ denote the $s$-panel containing $x_i$.
Now, by Property e) of the system of bijections, we have
$\beta(P,P') = \alpha(X_{n-1},X_n)\ldots \alpha(X_0,X_1) = \beta(P,P',sr_J)$,
where the last equality follows from Lemma \ref{le612nn}. As $w \in W_J$, it follows that
$w \prec sr_J$ and therefore the claim follows for this special case in view of Proposition~\ref{extension}.

We prove the claim for an arbitrary $w \in X_s$ by induction on $L_s(w)$ and
observe that the case $L_s(w) \leq 1$ is covered by the special case already considered before.
We put $t := w^{-1}sw$ and $Q:=\pi(P,w)=\pi(P',w)$ and remark that $Q$ is a $t$-panel.
 
%This is obvious if $w=1_W$, so assume $w \neq 1_W$.
%We will prove the result by induction on $L_s(w)$. Let $Q:=\pi(P,w)=\pi(P',w)$ be a $t$-panel (i.e. $w^{-1}sw=t$). By Lemma \ref{compatiblew}, %$L_s(P,Q)=L(P',Q)=L_s(w)$.

%Assume first that $L_s(w)=1$, so that $P$ and $Q$ are opposite in a rank 2 $J$-residue $R$. The panels $P'$ and
%$Q$ are also opposite in the same residue $R$. Hence $w=x_J=sr_J$. By the hypothesis on ${\cal B}$,
%$\aff(R)=A_f(R)=R\cap f^{\op}$ is connected. Choose $p\in P\cap f^{\op}$ and  $p'\in P'\cap f^{\op}$. There
%exists a gallery $\gamma$ from $p$ to $p'$ contained in $\aff(R)$. If $\gamma=(x_0=p,x_1,x_2,\ldots,x_n=p')$,
%let $X_i$ be the $s$-panel containing $x_i$. Those panels are in ${\cal P}^{\op}_s(f)$ and in $R$, and
%$X_i\sim_{t_i}X_{i+1}$ for some $t_i\in J$. By definition,
%$\beta(P,P')=\beta(\gamma,P,P'):=\alpha(X_{n-1},P')\ldots\alpha(X_1,X_2)\alpha(P,X_1)$. Moreover
%$\beta(X_i,X_{i+1})=\beta(X_i,X_{i+1},x_{J})$ if $t_i\neq s$ and  $\beta(X_i,X_{i+1})=1_{X_i}$ if $t_i=s$.
%Therefore $\beta(P,P')=\beta(P,P',x_J)$.

Assume $L_s(w)=k>1$ and assume that the result is proved for all $w'\in X_s$ with  $L_s(w')<k$. Let
$P=P_0,P_1,\ldots, P_k=Q$ be a compatible path from $P$ to $Q$, which exists by Lemma \ref{comppath}, and let
$P'=P'_0,P'_1,\ldots ,P'_k=Q$  be a compatible path from $P'$ to $Q$ with residues $R(P_i,P_{i+1})$ and
$R(P'_i,P'_{i+1})$ of the same type for all $0\leq i\leq k-1$. The existence of such a path follows from
Lemma \ref{compatiblew}.

Let $R$ be the rank 3 residue containing $R(P_{k-1},P_k)=R(P'_{k-1},P'_k)$ and
$R(P_{k-2},P_{k-1})$. Let $J$ be the type of $R$, which contains $t$ because $Q \subset R$. Let $T=\proj_{R}P$ and
$T'=\proj_{R'}P'$, which are panels by Lemma \ref{projpanels}. Let $c\in Q$. By Proposition \ref{pi}, $P\in
{\cal P}^{\op}_{s,c}(f)$, so there exists $x\in P\cap f_c^{\op}$. By Lemma \ref{projAfR}, $\proj_Rx\in A_f(R)$.
This means that $T$ contains chambers in $A_f(R)$, so whose $f$-value is $w_J$. We also have
$\delta(x,\proj_Rx) = w_J$ by Part d) of Proposition \ref{sa42} and as $\proj_R x \in \proj_R P = T$
we have $\proj_Tx = \proj_Rx$. It follows that
$\delta(x,\proj_{T}x)=w_J$, so that $T=\pi(P,w_J)$. By the same argument, $T'=\pi(P',w_J)$ and so $T'$ contains
chambers in $A_f(R)$. Therefore $T$ and $T'$ are $s'$-panels where $s'= w_J^{-1}sw_J\in J$.

Since $P_{k-2}$ is in a compatible path from $P$ (to $Q$), it is parallel to $P$. By Lemma \ref{projpanels}, $T$ is a
panel parallel to $P$ and $P_{k-2}$ and there exists a compatible path from $P$ to $P_{k-2}$ containing
$T$. Since all compatible paths between two given panels have the same length, the length of a compatible path
from $P$ to $T$ is less or equal to $k-2$. Hence $L_s(w_J)\leq L_s(w)-2$.

Choose $q\in T\cap A_f(R)$ and  $q'\in T'\cap A_f(R)$. Because of the hypothesis on ${\cal B}$, there exists a
gallery $\gamma$ from $q$ to $q'$ contained in $A_f(R)$. If $\gamma=(q=x_0,x_1,x_2,\ldots,x_n=q')$, let $X_i$
be the $s'$-panel containing $x_i$ and  $X_i\sim_{t_i}X_{i+1}$ for some $t_i\in J$ for $0\leq i\leq n-1$. By
Lemma~\ref{revpi}, there exists an $s$-panel $Q_i\in {\cal P}^{\op}_s(f)$ such that $X_i=\pi(Q_i,w_J)$ for all
$0\leq i\leq n$. Of course we take $Q_0=P$ and $Q_n=P'$. Since $\delta(Q_i,X_i)=w_J$, we have $\proj_RQ_i=X_i$
for all $0\leq i\leq n$. By Property $c)$ of Theorem \ref{beta},
$\beta(P,P')=\beta(Q_{n-1},Q_n)\ldots\beta(Q_1,Q_2)\beta(Q_0,Q_1)$.

Let $J_i=\{s',t_i\}\subset J$, let $R_i$ be the $J_i$-residue containing $X_i$ and $X_{i+1}$ and let
$w_i:=w_Js'r_{J_i} = sw^{J_i}$. 
%We have $w_i^{-1}sw_i=r_{J_i}s'w_J^{-1}sw_Js'r_{J_i}=r_{J_i}s'r_{J_i} \in J_i$, and
%$l(sw_Js'r_{J_i})=l(w_Jr_{J_i})=l(w_Js'r_{J_i})+1$ since $w_Jr_{J_i}=w_{J_i}r_{J_i}$ is the longest word of
%$w_JW_{J_i}$. Therefore $w_i\in X_s$ and $\pi(Q_i,w_i)=\pi(Q_{i+1},w_i)$ is the $t_i$-panel containing the only
%chamber of $R_i$ with $f$-value $w_{J_i}r_{J_i}$, that is $\proj_{R_i}f$. Therefore $Q_i\equiv_{w_i} Q_{i+1}$.
Now, $X_i = \pi(Q_i,w_J)$ and $X_{i+1} = \pi(Q_{i+1},w_J)$ are contained in the same spherical $J_i$-residue and
therefore $Q_i \equiv_{w_i} Q_{i+1}$ by Corollary \ref{ko611}.

Since $\proj_{R_i}Q_i=\proj_{R_i}\proj_RQ_i= \proj_{R_i}X_i=X_i$, a compatible path from $Q_i$ to $X_i$ (of
length $L_s(w_J)\leq k-2$) completed by the panel  $\pi(Q_i,w_i)$ is a compatible path of length $L_s(w_i)\leq
k-1$. By induction, this means that $\beta(Q_i,Q_{i+1})=\beta(Q_i,Q_{i+1},w_i)$.

Let $\tilde{w}:=w_Js'r_J=sw^J$. By Lemma \ref{le610nn} we have $\tilde{w} \in X_s$ and $w_i \prec \tilde{w}$.
It follows from Proposition~\ref{extension} that $Q_i \equiv_{\tilde{w}}Q_{i+1}$ and that 
$\beta(Q_i,Q_{i+1},w_i) = \beta(Q_i,Q_{i+1},\tilde{w})$. Hence
%Moreover, as $X_i =\pi(Q_i,w_J )$ and $X_{i+1} = \pi(Q_{i+1},w_J$ are contained in the $J$-residue $R$
%it follows that $Q 

%By a similar argument to the one for $w_i$, $\tilde{w}\in X_s$ and
%$\pi(Q_i,\tilde{w})$ is the $\tilde{w}^{-1}s\tilde{w}$-panel containing $\proj_{R}f$. Moreover
%$w_i\prec\tilde{w}$ for all $0\leq i\leq n-1$. Indeed
%$l(w_i^{-1}\tilde{w})=l(r_{J_i}s'w_J^{-1}w_Js'r_J)=l(r_{J_i}r_J)=l(r_J)-l(r_{J_i})=l(s'r_J)-l(s'r_{J_i})=l(w_J)+l(s'r_J)-(l(w_J)+l(s'r_{J_i}))=l(\tilde{w})-l(w_i)$.
%By Proposition \ref{extension}, we have  $Q_i\equiv_{\tilde{w}} Q_{i+1}$ and
%$\beta(Q_i,Q_{i+1},w_i)=\beta(Q_i,Q_{i+1},\tilde{w})$. Therefore
$$\begin{array}{ll}
\beta(P,P')&=\beta(Q_{n-1},Q_n)\ldots\beta(Q_1,Q_2)\beta(Q_0,Q_1)\\
&=\beta(Q_{n-1},Q_n,\tilde{w})\ldots\beta(Q_1,Q_2,\tilde{w})\beta(Q_0,Q_1,\tilde{w})\\
&= \beta(P,P',\tilde{w}).
\end{array}$$

We have $l(w^{-1}\tilde{w})=l(w^{-1}w_Js'r_J)=l(r_Js'w_J^{-1}w)=l(r_J)-l(s'w_J^{-1}w)$ because $s'w_J^{-1}w\in
W_J$. Moreover $l(s'w_J^{-1}w)=l(w_J^{-1}sw)=l(w_J^{-1}wt)=l(wt)-l(w_J)=l(w)+1-l(w_J)$ because $wt\in wW_J$.
Hence, on the one hand,  $l(w^{-1}\tilde{w})=l(r_J)+l(w_J)-1-l(w)$. On the other hand,
$l(\tilde{w})-l(w)=l(w_Js'r_J)-l(w)=l(w_J)+l(s'r_J)-l(w)=l(w_J)+l(r_J)-1-l(w)$ since $w_Js'r_J\in wW_J$.
Therefore $w\prec\tilde{w}$, and so  $\beta(P,P',w)=\beta(P,P',\tilde{w})$. This concludes the proof.
%By definition those panels are in ${\cal P}^{\op}_s(f)$ and in $R$, and $X_i\sim_{t_i}X_{i+1}$ for some $t_i\in J$.
%By definition, $\beta(P,P')=\beta(\gamma,P,P'):=\beta(X_{n-1},P')\ldots\beta(X_1,X_2)\beta(P,X_1)$.
%Moreover $\beta(X_i,X_{i+1})=\beta(X_i,X_{i+1},x_{J_i})$ where $J_i=\{s,t_i\}$.
% By Lemma \ref{Xs}, $x_{J_i}\prec x_J$, hence, by Proposition \ref{extension},  $X_i\equiv_{x_J} X_{i+1}$ and  $\beta(X_i,X_{i+1},x_{J_i})=\beta(X_i,X_{i+1},x_{J})$. Therefore $\beta(P,P')=\beta(P,P',x_J)$.
%On the other hand $w\prec x_J$, again by Lemma \ref{Xs}, and so  $\beta(P,P',w)=\beta(P,P',x_J)$. This allows us to conclude.
\end{Bw}

\begin{Ko}\label{caseII}
Let $R$ be a rank $2$ residue of ${\cal B}$, let $w_R$ be the shortest element
in $f(R)$ and suppose that $ w_R \in X_s$. Put $t:= w_R^{-1} s w_R \in J$. Let $c\in R$ and let $P,P'\in {\cal P}^{\op}_{s,c}(f)$. Then
$\beta(P,P')(\proj_Pc)=\proj_{P'}c$.
\end{Ko}
\begin{Bw}
Let $w=f(c)$, $J$ the type of $R$, so $w_R=w_J$. Let $d=\proj_Rf$, so that $f(d)=w^J$, where $w^J$ is the
unique longest word of $wW_J$. As  $w_J \in X_s$ and $w_J^{-1} s w_J \in J$  it follows that $sw^J\in X_s$
and
$u:= (sw^J)^{-1}s(sw^J) \in J$. Note that the $u$-panel $Q$ through $d$  is parallel to 
both $P$ and $P'$. Since $\delta(P,Q)=sw^J=w^Ju=\delta(P',Q)$ and $Q$ contains a chamber with $f$-value $w^Ju$, we
have $Q=\pi(P,sw^J)=\pi(P',sw^J)$. Therefore $P\equiv_{sw^J}P'$ and, by  Theorem \ref{betaw},
$\beta(P,P')=\beta(P,P',sw^J)=\proj_{P'}\proj_Q$.

Let $x\in P\cap f^{\op}_c$ and let $x'\in P'\cap f^{\op}_c$. By Lemma \ref{fopc}, there exist minimal
galleries $x=x_0,x_1,\ldots,x_n=c$ and $x'=x'_0,x'_1,\ldots,x'_n=c$, with $l(f(x_i))=i=l(f(x'_i))$ for all
$0\leq i\leq n$, containing $\proj_Pc$ and $\proj_{P'}c$, respectively. Obviously $\proj_Pc=x_1$ if $x_1\in P$ and
$x_0$ otherwise. By Proposition \ref{projf}, $f(c)=f(d)\delta(d,c)$; moreover
$l(f(c))=l(w)=l(w^J)-l(\delta(d,c))=l(f(d))-l(d,c)$. Hence there is a minimal gallery
$c=y_0,y_1,\ldots, y_m=d$ with $l(f(y_i))=l(f(c))+i$ for all $0\leq i\leq m$ and containing $\proj_Qc$,
where $m=l(c,d)$.
Obviously $\proj_Qc=y_{m-1}$ if $y_{m-1}\in Q$ and $y_m=d$ otherwise. We have that
$x=x_0,x_1,\ldots,x_n=y_0,y_1,\ldots, y_m=d$ is a minimal gallery, and so there is a minimal gallery (which is
a subgallery of the previous one) from $\proj_Pc$ to $\proj_Qc$ containing $c$. Therefore
$\proj_Pc=\proj_P\proj_Qc$. By a similar argument, $\proj_{P'}c=\proj_{P'}\proj_Qc$.

Putting everything together, $\beta(P,P')(\proj_Pc)=\proj_{P'}\proj_Q\proj_P\proj_Qc=\proj_{P'}\proj_Qc=\proj_{P'}c$,
because $P$ and $Q$ are parallel.
%By Lemma \ref{projAfR}, $\proj_Rx,\proj_Rx'\in A_f(R)$. Moreover $\delta(x,\proj_Rx)=w_J=\delta(x',\proj_Rx')$.
%As $l_f(R)\in X_s$, there is panels through $\proj_Rx$ and $\proj_Rx'$ of type $t$ and parallel to $P$, resp. $P'$. Because of $t=w_J^{-1}sw_J\in \typ(R)$, these panels are in $R$ and they are $\proj_RP=\pi(P,w_J)$ and $\proj_RP'=\pi(P',w_J)$.
\end{Bw}

\begin{Th}\label{main}
Let $c$ be a chamber of ${\cal B}$ and let $P,P'\in {\cal P}^{\op}_{s,c}(f)$. Then
$\beta(P,P')(\proj_Pc)=\proj_{P'}c$.
\end{Th}
\begin{Bw}
Throughout the proof we denote, for any residue $R$ of ${\cal B}$, by $w_R$ the 
unique shortest element
in the coset $f(R)$; this means that, if $R$ is a $J$-residue, and if
$w \in f(R)$, then $w_R = w_J$. Furthermore, we put $l_f(R):=l(w_R)=\min \{l(f(x)) \mid x \in R \}$.

We will prove the assertion by induction on $l(f(c))$.

Assume  $l(f(c))=0$. If $x\in f^{\op}_c$, then $\delta(x,c)=f(c)=1_W$, so  $f^{\op}_c=\{c\}$. Hence $P=P'$
contains $c$, and the statement is obvious since $\beta(P,P')=1_P$.

Assume  $l(f(c))=1$. If $f(c)=s$, then  $f^{\op}_c$ consists of all chambers $s$-adjacent to $c$ (except for
$c$ itself).  Hence $P=P'$ contains $c$, and the statement is again obvious. We now consider the case
$f(c)=t\neq s$. Let $R$ be the $\{s,t\}$-residue containing $c$. Then $w_R =1_W \in X_s$ and we are done by
Corollary \ref{caseII}.

Assume now  $l(f(c))=l\geq 2$ and assume that the theorem is proved for all chambers $c'$ with  $l(f(c'))<l$. Let $u,t$ be
the last two letters in a reduced word for $f(c)$, so that $l(f(c)ut)=l(f(c))-2$.  Let $R$ be the
$\{u,t\}$-residue containing $c$. Then $l_f(R)\leq l(f(c))-2$. 
If $\proj_RP$ is a panel, then $w_R \in X_s$ and we are done by Corollary \ref{caseII}.
Hence  we are left with the case where $\proj_RP$ is a chamber $p$ and 
$\proj_RP'$ is a chamber $p'$. Since $P$ contains a chamber $x$ in $f^{\op}_c$ and $\proj_RP=\proj_Rx$, we have
by Lemma \ref{projAfR} that $p\in A_f(R)$, and similarly $p'\in A_f(R)$. Moreover, there exists a minimal
gallery from $x$ to $c$ containing $p$
such that the length of the $f$-value strictly increases at each step, and so  $x\in f^{\op}_p$ by Lemma \ref{fopc}. Similarly $x'\in
f^{\op}_{p'}$. Recall that our general assumptions on ${\cal B}$
imply that the sets $A_f(R)$ are connected. Hence there exists a gallery $p=p_0,p_1,\ldots,p_n=p'$ (without
repetitions) entirely contained in $A_f(R)$. For all $1\leq j\leq n$, let $Q_j$ be the unique panel containing
$p_{j-1}$ and $p_j$, and let $z_j=\proj_{Q_j}f$. Since $l(f(z_j))=l_f(R)+1$, we have  $l(f(z_j))<l$ for all
$1\leq j\leq n$. For each  $1\leq j\leq n-1$, we can choose $x_j\in f^{\op}_{p_j}$.
We put $x_0:= x, x_n:=x'$ and denote the $s$-panel containing $x_i$ by $P_i$
for all $0 \leq i \leq n$.

By Lemma \ref{fopc}, there exists a gallery from $x_j$ to $p_j$ such that the length of the $f$-value strictly increases at each step,
for all  $1\leq j\leq n-1$. Since $l(f(z_j))=l(f(p_j))+1=l(f(z_{j+1}))$, $z_j,p_j\in Q_j$ and  $z_{j+1},p_j\in
Q_{j+1}$, by adding a chamber at the end of the previous gallery, we get two minimal galleries from $x_j$ to
$z_j$ and from $x_j$ to $z_{j+1}$, both such that the length of the $f$-value strictly increases at each step. Hence  $x_j\in
f^{\op}_{z_j}$ and  $x_j\in f^{\op}_{z_{j+1}}$ for all  $1\leq j\leq n-1$. Therefore  $P_j\in  {\cal
P}^{\op}_{s,z_j}(f)$ and  $P_{j}\in  {\cal P}^{\op}_{s,z_{j+1}}(f)$ for all  $1\leq j\leq n-1$. For a similar
reason $P=P_0\in  {\cal P}^{\op}_{s,z_1}(f)$ and  $P'=P_n\in  {\cal P}^{\op}_{s,z_n}(f)$. We conclude that
$P_{i-1},P_i\in  {\cal P}^{\op}_{s,z_i}(f)$ for all  $1\leq i\leq n$. Hence, by induction,
$\beta(P_{i-1},P_i)(\proj_{P_{i-1}}z_i)=\proj_{P_i}z_i$ for all  $1\leq i\leq n$.

As $w_R$ is not in $X_s$, it follows that
$\proj_RP_i$ is a chamber  for all  $1\leq i\leq n$. Because of this and since projections of residues on one another are
parallel,  $\proj_{P_i}R$ is also a chamber, hence $\proj_{P_i}z_i=\proj_{P_i}R=\proj_{P_i}c$ and
$\proj_{P_{i-1}}z_i=\proj_{P_{i-1}}R=\proj_{P_{i-1}}c$. Therefore we have
$\beta(P_{i-1},P_i)(\proj_{P_{i-1}}c)=\proj_{P_i}c$.

By Theorem \ref{beta}, $\beta(P,P')=\beta(P',P_{n-1})\ldots\beta(P_1,P_2)\beta(P,P_1)$
and therefore
$\beta(P,P')\proj_{P}c=\proj_{P'}c$.
\end{Bw}

\section{Adjacent codistances}\label{adjacentco}

In this section we will need the following facts.

\begin{Le} \label{Ropnn}
Let ${\cal B}$ be a spherical building of type $(W,S)$.
For each residue $R$ let $R^{\op}$ denote the set of all residues in ${\cal B}$  opposite $R$.
If $R$ and $T$ are two residues with $R^{\op} = T^{\op}$, then $R = T$.
\end{Le}
\begin{Bw}
Let $J$ be the type $R$ and consider a residue $R^*$ opposite $R$. Let $r \in W$ be the longest element.
 of $(W,S)$.
Then $R^*$ has type $K:= rJr$. As $R^*$ is opposite $T$, we have that $T$ is of type $J$.

Let $c \in R$ and put $d := \proj_Tc$. 
Put $w := \delta(c,d)$ and choose a chamber $e$ in ${\cal B}$ such that $\delta(d,e) =w_1 := w^{-1}r$
and consider the $K$-residue $R'$ containing $e$.
 Note first that $\delta(c,e) = r$ and that
$R'$ is therefore opposite $R$ and hence also opposite $T$ by our assumption.
We put $d' := \proj_Te$.
Now there exist a minimal gallery from $c$ to $d'$ passing through $d$, a minimal gallery
from  $d$ to $e$ passing through $d'$, and a minimal gallery from $c$ to $e$ passing through
$d$. We conclude that there is a minimal gallery from $c$ to $e$ passing through $d$ and $d'$.
As $T$  is opposite $R'$ and $e \in R'$ it follows that $\delta(e,d') = rr_J$ and therefore
$\delta(d',c) = r_J \in W_J$ because $\delta(e,c) = r$. This means that the $J$-residue containing $d'$ 
contains also $c$ and hence $R= T$.
\end{Bw}

The previous lemma has a very simple proof when considering buildings as simplical complexes, as originally defined by Tits \cite{Ti74}. Indeed, the residues $R$ and $T$ (which are just simplices in this setting) are, by definition of a `simplicial building', contained in an apartment, in which every simplex has a unique opposite. 

\begin{Le} \label{projcodistnn}
Let $f$ be a codistance on a building  ${\cal B}$, let $R$ be a spherical residue
and let $T$ be a residue contained in $R$. Then $\proj_Tf = \proj_T \proj_Rf$.
\end{Le}
\begin{Bw}
Let $c := \proj_Tf$ and $d := \proj_Rf$. Then we can find a minimal gallery from $c$ to $d$
such that the length of the $f$-value of the chambers in that gallery (strictly) increases at each step.
Hence $l(c,d) = l(f(d)) - l(f(c))$.
If $\proj_T d \neq c$, then we would have a chamber $e$ in $T$ such that 
$l(d,e) < l(d,c)$ which implies $l(f(e)) > l(f(c))$, which yields a contradiction.
\end{Bw}
 
%Put $d':= \proj_{R'} d$ and $w_2 := \delta(d,d')$.
%There exists a minimal gallery from $c$ to $e$ passing through $d$ and a minimal
%gallery from $d$ to $e$ passing through $d'$ because $d' := \proj_{R'} d$ and $e \in T$.
%Hence we have a minimal gallery from $c$ to $e$ passing through $d$ and through $d'$
%and therefore $l(r) = l(r,d) + l(d,d') + l(d',e)$.

 %Note first that $\delta(c,e) = r$ and that
%$R'$ is therefore opposite to $R$ and hence also opposite to $T$ by our assumption.

\begin{De}
Two codistances $f$ and $g$ on ${\cal B}$ are called {\it $s$-adjacent} if ${\cal P}^{\op}_{s}(f)={\cal
P}^{\op}_{s}(g)$. We denote this by $f\sim_s g$.
\end{De}

\begin{Le}\label{adjacency}
Let $f,g$ be two codistances on a building ${\cal B}$. Let $R$ be a spherical $J$-residue in $f^{\op}$. Let
$s\in J$. If $f$ and $g$ are $s$-adjacent, then $\proj_Rf$ and $\proj_Rg$ are $r_Jsr_J$-adjacent in ${\cal B}$.
\end{Le}
\begin{Bw}
Suppose $f$ and $g$ are $s$-adjacent. Then ${\cal P}^{\op}_{s}(f)={\cal P}^{\op}_{s}(g)$, which means that the
$s$-panels of $R$ in $f^{\op}$ and in $g^{op}$ coincide. The $r_Jsr_J$-panel $P$ containing $d:=\proj_Rf$ is
opposite in $R$ to all $s$-panels of $R$ in  $f^{\op}$.
Similarly the $r_Jsr_J$-panel $P'$ containing $d' :=\proj_Rg$ is opposite in $R$ to all $s$-panel of $R$ in $g^{\op}$.
By Lemma \ref{Ropnn} it follows that $P = P'$ and therefore $d$ and $d'$ are $r_Jsr_J$-adjacent.
% Suppose there is another panel $P'$ opposite to those
%same $s$-panels. Of course $P'$ is also of type  $r_Jsr_J$. Let $d':=\proj_{P'}f$. Then there exists a minimal
%gallery from $d$ to $d'$ which can be extended to a minimal gallery from $d$ to a chamber $c$ opposite to $d$,
%that is a chamber in $f^{\op}$. The $s$-panel containing $c$ is in $f^{\op}$, and so should be opposite to both
%$P$ and $P'$, which is not possible. Similarly,  the panel of type $r_Jsr_J$ containing $\proj_Rg$ is the only
%panel of $R$ opposite in $R$ to all $s$-panels in $g^{\op}$.  Therefore these 2 panels of type $r_Jsr_J$
%coincide, and so  $\proj_Rf$ and $\proj_Rg$ are $r_Jsr_J$-adjacent.
\end{Bw}

\begin{Le}\label{adjop}
Let $f$ be a codistance on a building ${\cal B}$, and let $g$ be a codistance $s$-adjacent to $f$.
Let $R$ be a $J$-residue in $f^{\op}$ with $s\in J$. Then $R$ is in $g^{\op}$.
\end{Le}
\begin{Bw}
Since $R$ is in  $f^{\op}$, $R$ contains a chamber $x$ in $f^{\op}$. The $s$-panel containing $x$ is in
$f^{\op}$, and so by hypothesis, it is in $g^{\op}$. Since this panel is in $R$, it means $R$ is in $g^{\op}$.
\end{Bw}

\begin{Le}\label{unicity}
Let $f$ be a codistance on a $k$-spherical building ${\cal B}$ such that $f^{\op}$ is connected, and let $g$ be
a codistance $s$-adjacent to $f$. Let $R$ be a $J$-residue of rank $\leq k-1$ in $f^{\op}$, with $s\in J$. Then
$\proj_Rg$ determines $g$ uniquely.
\end{Le}
\begin{Bw}
Suppose that $g_1$ and $g_2$ are two codistances $s$-adjacent to $f$ with   $\proj_Rg_1=\proj_Rg_2$.
By hypothesis, ${\cal P}^{\op}_{s}(f)={\cal P}^{\op}_{s}(g_1)={\cal P}^{\op}_{s}(g_2)$. %Let $d:=\proj_Rg_1=\proj_Rg_2$.

We claim that $g_1^{\op}\subseteq g_2^{\op}$. Let $x\in g_1^{\op}$. The $J$-residue $R_x$ containing $x$
is in $g_1^{\op}$. By Lemma \ref{adjop}, $R_x$ is also in $f^{\op}$.
%The $s$-panel containing $x$ is in $g_1^{\op}$, and so by hypothesis, it is in $f^{\op}$. Since this panel is in $R$, it means $R$ is in $f^{\op}$.
Since $f^{\op}$ is connected, there is a gallery $x_0,x_1,\ldots,x_n$ in $f^{\op}$ with $x_0\in R$ and $x_n\in
R_x$. We will show by induction on $n$ that $\proj_{R_x}g_1=\proj_{R_x}g_2$. If $n=0$, then $R_x=R$ and we are done.
Assume that we have shown that for every $J$-residue at "distance" (in the sense described
above) at most $n-1$ of $R$, the projections of $g_1$ and $g_2$ coincide. Let $R'$ be the $J$-residue
containing $x_{n-1}$. By the induction hypothesis, $\proj_{R'}g_1=\proj_{R'}g_2$. We have  $x_{n-1}\sim_t x_n$.
If $t\in J$ then $R_x=R'$ and we are done. So assume $t\not\in J$, let $K=J\cup\{t\}$ (which is spherical by
hypothesis) and let $\tilde{R}$ be the $K$-residue containing $R_x$ and $R'$. Then, by Lemma \ref{Ropnn}, the residue of type
$\op_K(J):=\{r_{K}ur_K|u\in J\}$ containing $\proj_ {\tilde{R}}f$ is the unique  residue of  $\tilde{R}$ opposite
in  $\tilde{R}$ to all $J$-residues of  $\tilde{R}$ in  $f^{\op}$.
Similarly,  the $\op_K(J)$-residue containing $\proj_{\tilde{R}}g_i$ is the only residue of
${\tilde{R}}$ opposite in ${\tilde{R}}$ to all $J$-residues in $g_i^{\op}$, for $i=1,2$. By Lemma \ref{adjop},
the sets of $J$-residues in $f^{\op}$ and in $g_i^{\op}$ ($i=1,2$) coincide. Therefore these three
$\op_K(J)$-residues coincide;  let us name it $T$.
%Then by Lemma \ref{adjacency},  $\proj_Rf$ and $\proj_Rg_1$ are $r_Jsr_J$-adjacent in ${\cal B}$, as well as  $\proj_Rf$ and $\proj_Rg_2$. Let $T$ be the $r_Jsr_J$-panel containg these 3 chambers, which is parallel to $P'$ and to $Q$.
As $f(y)=f(\proj_{\tilde{R}}f)\delta(\proj_{\tilde{R}}f,y)$ for $y\in \tilde{R}$, by Lemma \ref{projf}, we have
$l(f(y))=r_K-l(\proj_{\tilde{R}}f,y)$ for $y\in\tilde{R}$, and so
$\proj_{R'}f=\proj_{R'}\proj_{\tilde{R}}f$. Similarly for $g_1,g_2$ and for $R_x$. We have
$\proj_{R'}g_1=\proj_{R'}g_2$, and so $\proj_{R'}\proj_{\tilde{R}}g_1=\proj_{R'}\proj_{\tilde{R}}g_2$, with
$\proj_{\tilde{R}}g_1,\proj_{\tilde{R}}g_2\in T$. 
As $T$ contains $\proj_{\tilde{R}} f, \proj_{\tilde{R}} g_i (i=1,2)$ 
and since
$R'$ and $T$ are parallel we have
$\proj_{\tilde{R}}g_1=\proj_{\tilde{R}}g_2$. Now  $\proj_{R_x}g_1 = \proj_{R_x}\proj_{\tilde{R}} g_1 =
\proj_{R_x}\proj_{\tilde{R}} g_2 = \proj_{R_x}g_2$
Since $x\in g_1^{\op}$ and, for any $y\in R_x$,
$g_1(y)=g_1(\proj_{R_x}g_1)\delta(\proj_{R_x}g_1,y)$  by Lemma \ref{projf}, we have
$1_W=r_J\delta(\proj_{R_x}g_1,x)$. Therefore $r_J=\delta(\proj_{R_x}g_1,x)=\delta(\proj_{R_x}g_2,x)$, which
implies that  $x\in g_2^{\op}$. By symmetry, we get  $g_1^{\op}=g_2^{\op}$. We now conclude by Lemma
\ref{fopunique}.
\end{Bw}

\begin{Sa} Let $\tilde{{\cal C}}$ be the set of all codistances on a 3-spherical building ${\cal B}$.
Then $(\tilde{{\cal C}},(\sim_s)_{s\in S})$ is a chamber system.
\end{Sa}
\begin{Bw}
It follows from the definition that $\sim_s$ is an equivalence relation on  $\tilde{{\cal C}}$ for all $s\in
S$. Suppose $f\sim_s g$ and $f\sim_t g$ for $s,t\in S$ and $f\neq g$. Let $J=\{s,t\}$, which is spherical. Let
$R$ be a $J$-residue in $f^{\op}$. By Lemma \ref{unicity}, $\proj_Rf$ and $\proj_Rg$ are distinct. By Lemma
\ref{adjacency}, the chambers $\proj_Rf$ and $\proj_Rg$ are $r_Jsr_J$-adjacent and also $r_Jtr_J$-adjacent in
${\cal B}$. Since the chambers of ${\cal B}$ form a chamber system, it means $r_Jsr_J=r_Jtr_J$, and hence
$s=t$.
\end{Bw}

From now on, we again assume that ${\cal B} =({\cal C},\delta)$ is a 3-spherical building of type $(W,S)$
satisfying (lco) and (lsco) and that $f$ is a codistance on ${\cal B}$. Let ${\cal B}^*=({\cal
C}^*,(\sim_s)_{s\in S})$ be the chamber system on the connected component of $f$.

Fix $s\in S$ and $\tilde{P}$ in  ${\cal P}^{\op}_{s}(f)$. For each chamber $p$ of $\tilde{P}$ in $f^{\op}$, we
will define another codistance on ${\cal B}$. Let $\beta(p):=\{\beta(\tilde{P},Q)(p)|Q\in {\cal
P}^{\op}_{s}(f)\}$. By Theorem \ref{beta}, this set contains exactly one chamber in each panel of  ${\cal
P}^{\op}_{s}(f)$, none of which is the projection of $f$ on it.

\begin{Th}\label{existence}
For $c\in {\cal B}$, choose  $P\in{\cal P}^{\op}_{s,c}(f)$, and put
$$g(c)=
\left\{\begin{array}{ll}
sf(c) &\text{ if } \proj_Pc\in\{\proj_Pf,\beta(p)\cap P\}\\
f(c) & \text{ otherwise. }
\end{array}\right.$$
Then $g$ is a codistance on ${\cal B}$. Moreover $g$ is $s$-adjacent to $f$ and, for  $P\in{\cal
P}^{\op}_{s}(f)$, $\proj_Pg=\beta(p)\cap P$.
\end{Th}
\begin{Bw}
The function $g:{\cal C}\rightarrow W$ is independent of the choice of $P$ by Theorem \ref{main} and by statement d) of Theorem
\ref{beta}. Let $Q$ be a $t$-panel, so that  $f(x)\in \{w,wt\}$ for all $x\in Q$ and $Q$ contains a unique
chamber $q:=\proj_Qf$ with $f$-value the longest word of the two, which we can assume to be $wt$. We first show that $Q$ satisfies the codistance condition for $f$. We distinguish two cases. 

Case 1. Assume $w^{-1}sw=t$. Then $w\in X_s$ and there exists $P\in{\cal P}^{\op}_{s}(f)$ parallel to $Q$ with
$\delta(P,Q)=w$. By Proposition~\ref{pi},   $P\in{\cal P}^{\op}_{s,x}(f)$ for all chambers $x$ of $Q$. Since
$P$ and $Q$ are parallel, $\proj_P^Q$ and $\proj_Q^P$  are inverse bijections between $P$ and $Q$. Hence
$g(x)=f(x)$ for all $x\in Q$, except for $q$ whose $g$-value is $sf(q)=swt=w$ and for $\proj_Q(\beta(p)\cap P)$
whose $g$-value is $sf(\proj_Q(\beta(p)\cap P))=sw=wt$. Hence  $g(x)\in \{w,wt\}$  and $Q$ contains a unique
chamber with $g$-value $wt$.

Case 2. Now assume  $w^{-1}sw\neq t$. We claim that either $g(x)=f(x)$ for all $x\in Q$ or $g(x)=sf(x)$ for all $x\in
Q$. Suppose we proved the claim. In the first case, it is obvious that $Q$ will satisfy the codistance condition for $g$. Suppose we are in the
second case. Then $g(x)\in \{sw,swt\}$ for all $x\in Q$ and $Q$ contains a unique chamber with $g$-value $swt$.
We just need to show that $l(swt)=l(sw)+1$ to conclude that $Q$ satisfies the codistance condition for $g$. If
$l(sw)=l(w)+1$, it follows from Assertion b) of Proposition
\ref{sa41}  that either $l(swt)=l(w)+2$ or $swt=w$. Since the second case is excluded, we have
$l(swt)=l(w)+2=l(sw)+1$. If $l(sw)=l(w)-1$ and $l(swt)=l(sw)-1$, then $l(swt)=l(w)-2=l(wt)-3$, and we get a
contradiction, hence if $l(sw)=l(w)-1$ we also get $l(swt)=l(sw)+1$.

We now prove the claim. Let $x\in Q$ with $f$-value $w$, $y\in f^{\op}_x$ and $P$ the $s$-panel containing $y$,
so that  $P\in{\cal P}^{\op}_{s,x}(f)$. If we add the chamber $q$ to a minimal gallery from $y$ to $x$, we get
a minimal gallery from $y$ to $q$ with the required condition on $f$, and so, by Lemma \ref{fopc},  $y\in
f^{\op}_q$ and $P\in{\cal P}^{\op}_{s,q}(f)$. Let $x'$ be another chamber of $Q$ with $f$-value $w$ and
$P'\in{\cal P}^{\op}_{s,x'}(f)$. By the same argument, $P'\in{\cal P}^{\op}_{s,q}(f)$. By Theorem \ref{main},
this means $\beta(P,P')(\proj_Pq)=\proj_{P'}q$.  Since $P$ and $Q$ ($P'$ and $Q$, respectively) are not parallel,
$\proj_PQ$ ($\proj_{P'}Q$, respectively) is a chambers, and so $\proj_Pq=\proj_Px$ ($\proj_{P'}q=\proj_{P'}x'$, respectively).
Therefore  $\beta(P,P')(\proj_Px)=\proj_{P'}x'$, and so $\proj_Px\in\{\proj_Pf,\beta(p)\cap P\}$ if and only if
$\proj_{P'}x'\in\{\proj_{P'}f,\beta(p)\cap P'\}$ by Theorem \ref{beta} d). 
Moreover we also have  $\proj_Px\in\{\proj_Pf,\beta(p)\cap
P\}$ if and only if $\proj_{P}q\in\{\proj_Pf,\beta(p)\cap P\}$. Therefore the claim is proved. 

Hence we have shown that $g$ is a codistance. We now show that $g$ is $s$-adjacent to $f$. 

Let $P\in{\cal P}^{\op}_{s}(f)$. Then $P$ contains chambers in $f^{\op}$ and one chamber $p$ with $f(p)=s$.
Obviously $P\in {\cal P}^{\op}_{s,p}(f)$, hence $\proj_Pp=\proj_Pf = p$ and so $g(p)=sf(p)=1_W$. Hence $p\in g^{\op}$ and
$P\in{\cal P}^{\op}_{s}(g)$. Let $P$ be a $s$-panel not in ${\cal P}^{\op}_{s}(f)$. Then $f(x)\in\{w,ws\}$ for
$x\in P$ with $s\neq w\neq 1$. Hence $g(x)\in\{w,sw,ws,sws\}$ for $x\in P$. Since $1_W\not\in\{w,sw,ws,sws\}$,
no chamber of $P$ is in $g^{\op}$, and so  $P\not\in{\cal P}^{\op}_{s}(g)$. This proves that $f\sim_s g$.

Finally, let $P\in{\cal P}^{\op}_{s}(f)$. Then for any $c\in P$, $P\in{\cal P}^{\op}_{s,c}(f)$, therefore
$g(c)=f(c)$ unless $c\in \{\proj_Pf,\beta(p)\cap P\}$. Hence the only chamber of $P$ with $g$-value $s$ is $\beta(p)\cap P$.
\end{Bw}

\begin{Sa}\label{alpha} Let ${\cal B}$ be a 3-spherical building of type $(W,S)$ satisfying {\rm (lco)} and {\rm (lsco)}.
Let $J\subseteq S$ be spherical, and let $f$ be a codistance on  ${\cal B}$.
Let $R$ be a $J$-residue of ${\cal B}$ in $f^{\op}$ and let $\tilde{R}$ be the $J$-residue containing $f$ in  ${\cal B}^*$.
Then $\alpha:\tilde{R}\rightarrow R:g\rightarrow \proj_Rg$ is a bijection such that:\\
$(i)$ $\forall g_1,g_2\in \tilde{R},s\in J$, we have $g_1\sim_s g_2$ if and only if
$\alpha(g_1)\sim_{r_Jsr_J}\alpha(g_2)$, \\
$(ii)$  $\forall g\in \tilde{R},c\in R$, we have $c\in g^{\op}$  if and only if $\delta(\alpha(g),c)=r_J$.
\end{Sa}
\begin{Bw}
By Lemma \ref{unicity}, $\alpha$ is injective. Let $d=\proj_Rf$. We will show by induction on the numerical distance
$l(x,d)$ that $x \in \alpha(\tilde{R})$. First notice  that $\alpha(f)=d$, so  $x \in \alpha(\tilde{R})$  if
$l(x,d)=0$. Suppose we have proved that  $x \in \alpha(\tilde{R})$  for all $x$
satisfying $l(x,d)<l$ and suppose $(ly,d)=l$. Let $y=y_0,y_1,\ldots,y_l=d$ be a minimal
gallery. By hypothesis, there exists $g_1 \in \tilde{R}$ with $\alpha(g_1)=y_1$. Let $T$ be the $t$-panel
containing $y_0$ and $y_1$ for some $t\in S$. Let $s=\op_J(t)=r_Jtr_J\in J$. 
By Lemma~\ref{adjop},
$R$ is in $g^{\op}$ for any $g\in \tilde{R}$ and so in particular for $g_1$. Therefore there exists
$c\in R\cap g_1^{\op}$ and the $s$-panel  $P$ containing $c$ is in  ${\cal P}^{\op}_{s}(g_1)$. By construction
$P$ and $T$ are opposite and hence parallel. Let $p:=\proj_Py$. Using Theorem \ref{existence}, we can construct a
codistance $g$ which is $s$-adjacent to $g_1$ with $\proj_Pg=p$. By Lemma~\ref{adjacency}, $\proj_Rg$ and $\proj_Rg_1$
are $t$-adjacent, and so $\proj_Rg\in T$. Since $\proj_Pg=\proj_P\proj_Rg$ by
Lemma~\ref{projcodistnn}, we must
have $\proj_Rg=y$. Therefore $\alpha(g)=y$ and $\alpha$ is surjective.

By Lemma \ref{adjacency}, if $g_1$ and $g_2$ are $s$-adjacent in $\tilde{R}$, then  $\proj_Rg_1$ and
$\proj_Rg_2$ are $r_Jsr_J$-adjacent in $R$.

Now assume  $g_1$ and $g_2$ are codistances in $\tilde{R}$ with  $\proj_Rg_1\sim_{r_Jsr_J}\proj_Rg_2$ for some
$s\in J$. Let $P$ be the $r_Jsr_J$-panel containing them and put $e:=\proj_Pd$. As $\alpha$ is surjective, there
exists $g\in\tilde{R}$ with $\alpha(g)=e$. We have shown above that there exist codistances $g_1'$ and $g_2'$,
both $s$-adjacent to $g$, with $\proj_Pg_1'=\proj_Pg_1$ and $\proj_Pg_2'=\proj_Pg_2$. By the injectivity of
$\alpha$, $g_1'=g_1$ and $g'_2=g_2$, and so $g_1$ and $g_2$ are both  $s$-adjacent to $g$. Since ${\cal B}^*$
is a chamber system, this means  $g_1\sim_sg_2$. This proves $(i)$.

We now prove $(ii)$. Let  $g\in \tilde{R}$. By Lemma \ref{projf}, for all $c\in R$,
$g(c)=g(\alpha(g))\delta(\alpha(g),c)$. Since $R\in g^{\op}$ as noticed above, $g$ takes on $R$ its values in
$W_J$, and so  $g(\alpha(g))=r_J$. Hence $c\in g^{\op}$ if and only if $g(c)=1_W$ if and only if $\delta(\alpha(g),c)=r_J$.
\end{Bw}

\begin{Ko}\label{diagram} The chamber system ${\cal B}^*$ has the same diagram as  ${\cal B}$.
\end{Ko}
\begin{proof}
Let $M$ be the diagram of ${\cal B}$, which mean that each rank 2 $J$-residue is a generalized $M_J$-gon.
Let $\tilde{R}$ be a $J$-residue of rank 2 of ${\cal B}^*$. %, which is spherical by hypothesis.
Let $g$ be a codistance in $\tilde{R}$ and let $R$ be a $J$-residue in $g^{\op}$. Then, by Proposition
\ref{alpha}, $\tilde{R}$ is a building of the same type as $R$, hence a generalized $M_J$-gon. Therefore ${\cal
B}^*$ has diagram $M$.
\end{proof}

\section{Construction of the twinning}\label{constructiontwinning}

In order to construct a twinning we apply the main result
of \cite{Mu98} which we recall below and whose statement
requires some preparation.

Let $(W,S)$ be a Coxeter system and let ${\cal B}_+ = ({\cal C}_+,\delta_+),
{\cal B}_- = ({\cal C}_-,\delta_-)$ be two buildings of type $(W,S)$.
%A \emph{twinning} of ${\cal B}_+$ and ${\cal B}_-$ is a mapping $\delta_*$ blablabla.
An \emph{opposition relation} between ${\cal B}_+$ and ${\cal B}_-$ is a non-empty subset ${\cal O}$ of ${\cal
C}_+ \times {\cal C}_-$ such that there exists a twinning $\delta_*$ of ${\cal B}_+$ and ${\cal B}_-$ with the
property that ${\cal O} = \{ (x,y) \in {\cal C}_+ \times {\cal C}_- \mid \delta_*(x,y) = 1_W \}$.

A \emph{local opposition relation}  between ${\cal B}_+$ and ${\cal B}_-$ is
a non-empty subset ${\cal O}$ of ${\cal C}_+ \times {\cal C}_-$
such that for each $(x,y) \in {\cal O}$ and each subset
$J\subseteq S$ of cardinality at most 2 the set ${\cal O} \cap (R_J(x) \times R_J(y))$
is an opposition relation between the $J$-residues of $x$ and $y$.
Note that the definition of a local opposition relation
makes perfect sense for two chamber systems of type $(W,S)$
as well.

Here is the main result of \cite{Mu98}.
\begin{Th}
Let $(W,S)$ be a Coxeter system and let ${\cal B}_+ = ({\cal C}_+,\delta_+),
{\cal B}_- = ({\cal C}_-,\delta_-)$ be two thick buildings of type $(W,S)$
and let ${\cal O}$ be a non-empty subset
of ${\cal C}_+ \times {\cal C}_-$. Then ${\cal O}$ is an opposition relation
between ${\cal B}_+$ and ${\cal B}_-$ if and only if it is a local
opposition relation between the two buildings.
\end{Th}

The following corollary of the previous theorem has been proved
in \cite[p.28]{Mu99}. We paraphrase that proof here.

\begin{Ko}\label{HabilKo} Let $(W,S)$ be a Coxeter system and let
$({\cal C}_+,(\sim_s)_{s \in S}), ({\cal C}_-,(\sim_s)_{s \in S})$
be two connected, thick chamber systems of type $(W,S)$ whose
universal $2$-covers are buildings. Suppose that there exists
a local opposition relation ${\cal O} \subseteq ({\cal C}_+ \times {\cal C}_-)$
between them. Then the chamber systems are buildings.
In particular, there exist unique distances
$\delta_+: {\cal C}_+  \times {\cal C}_+  \rightarrow W$,
$\delta_-: {\cal C}_-  \times {\cal C}_- \rightarrow W$ and
$\delta_*:  ({\cal C}_+  \times {\cal C}_-) \cup
({\cal C}_- \times {\cal C}_+)  \rightarrow W$ such that
$(({\cal C}_+,\delta_+), ({\cal C}_-,\delta_-),\delta_*)$
is a twin building of type $(W,S)$ and such that
${\cal O} = \{ (x,y) \in {\cal C}_+  \times {\cal C}_- \mid
\delta_*(x,y) = 1_W \}$.
\end{Ko}
\begin{Bw}
Let $\overline{{\cal B}}_{\epsilon}=(\overline{{\cal C}}_{\epsilon},(\sim_s)_{s \in S})$ be the universal
2-cover of $({\cal C}_{\epsilon},(\sim_s)_{s \in S})$, which is a building by hypothesis, with covering
morphism $\phi_{\epsilon}:\overline{{\cal C}}_{\epsilon}\rightarrow {\cal C}_{\epsilon}$, for $\epsilon=+,-$.
Let $\overline{{\cal O}}=\{(x,y)\in \overline{{\cal C}}_+\times\overline{{\cal C}}_-|(\phi_+(x),\phi_-(y)\in
{\cal O} \}$. Obviously $\overline{{\cal O}}$ is a local opposition relation between $\overline{{\cal B}}_+$
and $\overline{{\cal B}}_-$. By the previous theorem, this means that  $\overline{{\cal O}}$ is the opposition
relation of a twin building $(\overline{{\cal B}}_+,\overline{{\cal B}}_-,\overline{\delta}_*)$.

Let $\overline{x}\neq\overline{y}\in \overline{{\cal C}}_-$. By Lemma~\ref{+*}, $\overline{x}^{\op}\neq
\overline{y}^{\op}$. Hence there exists $\overline{z}\in \overline{{\cal C}}_+$ such that
$(\overline{z},\overline{x})\in\overline{{\cal O}}$ but $(\overline{z},\overline{y})\not\in\overline{{\cal
O}}$. If $\phi_-(\overline{x})=v=\phi_-(\overline{y})$, then we have both $(\phi_+(\overline{z}),v)\in {\cal
O}$ and  $(\phi_+(\overline{z}),v)\not\in {\cal O}$, a contradiction. This shows that $\phi_-$ is injective and
hence is the identity. The same argument shows that $\phi_+$ is the identity. Therefore $\overline{{\cal
C}}_{\epsilon}={\cal C}_{\epsilon}$ for $\epsilon=+,-$ and the result follows.
\end{Bw}

In order to apply the corollary above we need the following
lemma.

\begin{Le} \label{le73}
Let ${\cal B}_+ = ({\cal C}_+,\delta_+),
{\cal B}_- = ({\cal C}_-,\delta_-)$ be two buildings of spherical
type $(W,S)$, let $r \in W$ be the longest element in $W$
and let ${\cal O}$ be a non-empty subset
of ${\cal C}_+ \times {\cal C}_-$. Then the following are
equivalent.
\begin{itemize}
\item[a)] ${\cal O}$ is an opposition relation between
${\cal B}_+$ and ${\cal B}_-$.
\item[b)] There exists a bijection $\alpha:{\cal C}_+ \rightarrow {\cal C}_-$
such that the following two conditions are satisfied:
\begin{itemize}
\item[(i)] For all $x,y \in {\cal C}_+$ and all $s \in S$ we have $x \sim_s y$
if and only if $\alpha(x) \sim_{rsr} \alpha(y)$;
\item[(ii)] ${\cal O} = \{ (x,y) \in {\cal C}_+ \times {\cal C}_- \mid
\delta_-(\alpha(x),y) = r \}$.
\end{itemize}
\end{itemize}
\end{Le}
\begin{Bw}
Suppose ${\cal O}$ is an opposition relation between ${\cal B}_+$ and ${\cal B}_-$. Then there exists a
twinning $\delta_*: {\cal C}_+ \times {\cal C}_- \rightarrow W$ inducing the opposition relation ${\cal O}$. 
Let $x$ be a chamber in ${\cal C }_+$ and let $f_x: {\cal C}_- \rightarrow W$ be defined by
$f_x(y) := \delta_*(x,y)$. Then $f$ is a codistance on ${\cal B}_-$  and as ${\cal B}_-$ is
spherical, $\proj_{{\cal C}_-} f_x$ makes sense. It is the unique chamber in ${\cal C}_-$ at codistance
$r$ to $x$, where $r$ denotes the longest element in $W$.
 Define  $\alpha:{\cal C}_+ \rightarrow {\cal C}_-$ by $\alpha(x)=\proj_{{\cal C}_-}f_x$.
One checks that $\alpha$ is a
bijection and satisfies $(i)$ and $(ii)$.

Now suppose there exists a bijection $\alpha$ satisfying $(i)$ and $(ii)$. We define a mapping $\delta_*$ from
$({\cal C}_+ \times {\cal C}_-)\cup ({\cal C}_- \times {\cal C}_+)$ into $W$ by
$\delta(x,y):=r\delta_-(\alpha(x),y)$ and  $\delta(y,x):=\delta(x,y)^{-1}$, for $x\in{\cal C}_+$ and $y
\in{\cal C}_-$. Using the axioms of buildings, it can easily be checked that $\delta_*$ is a twinning. Moreover
${\cal O} = \{ (x,y) \in {\cal C}_+ \times {\cal C}_-\mid \delta_*(x,y) = 1_W \}$, so ${\cal O}$ is an
opposition relation between ${\cal B}_+$ and ${\cal B}_-$.

\end{Bw}

\subsection*{End of the proof of the main result}
Let ${\cal B}_- = ( {\cal C}_-,\delta_-)$ be a thick building
of type $(W,S)$ satisfying
all necessary properties and let $f:  {\cal C}_- \rightarrow W$ be
a codistance.

Consider the chamber system of all codistances of ${\cal B}_-$
which is a chamber system over $S$.
Let ${\cal C}_+$ be the connected component containing $f$ and
consider the chamber system $({\cal C}_+, (\sim_s)_{s \in S})$ which
is a connected chamber system of type $(W,S)$
by Corollary \ref{diagram}. It readily follows from Proposition
\ref{alpha} that all $J$-residues of rank at most 3 are spherical
buildings and in particular that
$({\cal C}_+, (\sim_s)_{s \in S})$ is thick .
By a result of Tits \cite{Ti81}  it
follows that the universal 2-cover of this chamber system is a building.

We define ${\cal O} \subseteq {\cal C}_+ \times {\cal C}_-$
by setting ${\cal O} := \{ (g,c) \in  {\cal C}_+ \times {\cal C}_- \mid
g(c) = 1_W \}$. Using Lemma \ref{le73} and Proposition \ref{alpha}
we see that ${\cal O}$ is a local opposition between
the  chamber systems $({\cal C}_+, (\sim_s)_{s \in S})$
and $({\cal C}_-, (\sim_s)_{s \in S})$, which are both thick chamber systems
of type $(W,S)$ whose universal covers are buildings. Therefore,
Corollary \ref{HabilKo} yields the twin building.

Now we have $f':=\delta_*(f,.)$ is a codistance on ${\cal B}_-$ with
$$f'^{\op}=\{c\in {\cal C}_-|\delta_*(f,c)=1_W\}=\{c\in {\cal C}_-|(f,c)\in
{\cal O}\}=\{c\in {\cal C}_-|f(c)=1_W\}=f^{\op}.$$
By Lemma \ref{fopunique}, we have $f'=f$ and so  $\delta_*(f,x)=f(x)$ for all $x\in {\cal C}_-$.

\section{Remarks on the conditions in the main result}

The purpose of this  section is to  provide some additional information about 
the conditions on the buildings in our main result. In the discussion below we always
assume the buildings are of irreducible type which is not a serious restriction,
because the general case can be reduced to the irreducible case.

%\smallskip
%\noindent
%{\bf 3-sphericity:} It is fairly easy to construct examples of buildings admitting a codistance which cannot be realized
%as a `half of a twin building'. For instance, it is a trivial fact that 
%each thick building  ${\cal B}_-$ of type $\tilde{A}_1$ admits a codistance 
%$f$. Moreover, 
%it can be shown that the conclusion of the main result holds 
%for a thick  $\tilde{A}_1$-building if and only if
%panels of the same type have the same cardinality (see \cite{AB99}, \cite{RT99}).
%It is most likely, that this remains true for all right-angled buildings.
%
%If there are finite entries different from 2 in the diagram, the question becomes more delicate. Nevertheless,
%we expect a similar behavior if there are `enough' infinities in the diagram. So it is natural
%to assume that the diagram is {\sl 2-spherical} which means that there are no infinities in the diagram.
%By the following remarks, the conditions asked in addition to 3-sphericity are `almost always'
%satisfied and therefore it remains to consider 2-spherical buildings which are not 3-spherical.
%We have no idea about what to expect in this case. On the one hand, the methods used in  the proof
%of our main result completely fail in this more general context. On the other hand we could not manage
%to construct counter-examples in the $\tilde{A}_2$-case - a case which is well understood in a lot of respects. 

\smallskip
\noindent
{\bf 3-sphericity:} % (proposed by Alice) 
If we drop the 3-sphericity condition (together with conditions (lco) and (lsco)), the conclusion of our main result 
is not always true. Indeed it is fairly easy to construct examples of buildings admitting a codistance which cannot be realized as a `half of a twin building'. For instance, it is a trivial fact that each thick building  ${\cal B}_-$ of type $\tilde{A}_1$ admits a codistance 
$f$. Moreover, 
it can be shown that  ${\cal B}_-$ can be realized as a `half of a twin building' if and only if
panels of the same type have the same cardinality (see \cite{AB99}, \cite{RT99}).

It is an interesting question to wonder which buildings admitting a codistance can or cannot be realized as a `half of a twin building'.
It is most likely that all right-angled buildings admit a codistance, and that they  can be realized as a `half of a twin building' if and only if
panels of the same type have the same cardinality. If there are finite entries different from 2 in the diagram, the question becomes more delicate. Nevertheless, we expect a behavior  similar to the case of right-angled buildings if there are `enough' infinities in the diagram.
Hence, for the conclusion of our main result to hold,  it is natural
to assume that the diagram is {\sl 2-spherical} (i.e. there are no infinities in the diagram), in which case panels of the same type always have the same cardinality.
By the following remarks, the conditions asked in addition to 3-sphericity are `almost always'
satisfied and therefore it remains to consider 2-spherical buildings which are not 3-spherical.
We have no idea about what to expect in this case. On the one hand, the methods used in  the proof
of our main result completely fail in this more general context. On the other hand we could not manage
to construct counter-examples in the $\tilde{A}_2$-case --- a case which is well understood in a lot of respects. 

\smallskip
\noindent
{\bf Condition (lco):} By the 3-sphericity assumption, all entries in the diagram
are equal to 2, 3 or 4, if the rank is at least 3.  
It follows from an observation of Cuypers, see \cite{Br93}, that Condition
(lco) is satisfied if there is no rank 2 residue isomorphic
to the building associated with $B_2(2)$. In particular, Condition (lco)
is satisfied if the diagram is simply laced. 

\smallskip 
\noindent
{\bf Condition (lsco):} 
 It follows from \cite{Ti86} Corollaire 2 that Condition (lsco)
is satisfied if the diagram is simply laced and if each panel contains at
least 4 chambers. If there are subdiagrams of type $B_2$ we have to
consider buildings of type $B_3$. For those
the relevant results concerning Condition (lsco) may be
found in  \cite{Ab96}. They imply that Condition (lsco)
is satisfied if each residue of type $B_3$ comes from an embeddable
polar space and if each panel contains at least 17 chambers.
The first condition is equivalent to the fact that any $A_2$-residue
corresponds to a Desarguesian projective plane, and it is
very likely that it can be dropped. Moreover, it is expected
that the bound 17 is not optimal.

\end{document}